\documentclass[10 pt, leqno]{amsart}

\usepackage {amsfonts}
\usepackage{amsthm}
\usepackage{amssymb}
\usepackage{latexsym}
\usepackage{amsmath}
\usepackage{mathrsfs}

\theoremstyle{plain}
\newtheorem{tw}{Theorem}[section]

\newtheorem {lem} [tw]{Lemma}

\newtheorem{cor}[tw]{Corollary}

\theoremstyle{definition}
\newtheorem {deft}[tw] {Definition}
\newtheorem {rem} [tw]{Remark}

\newcommand{\Pois}{R_{\Vcont}}
\newcommand{\sPois}{R_{s,\Vcont}}
\newcommand{\cPois}{T_{\Vcont}}

\newcommand{\bc} {\mathbb C}

\newcommand{\bn}{\mathbb N}
\newcommand{\br}{\Bbb R}

\newcommand{\bt}{\Bbb T}

\newcommand{\alg} {\mathsf{A}}

\newcommand {\Pop}{\Delta_s(\Vcont)}
\newcommand {\Wop}{W_s(\Vcont)}

\newcommand {\Fix} {{\textup{Fix}}\;}

\newcommand{\Hil}{\mathsf{H}}

\newcommand{\Kil}{\mathsf{K}}

\newcommand{\Vcont}{\mathcal{V}}
\newcommand{\cont}{\mathcal{C}_{\La}}
\newcommand{\Wcont}{\mathcal{W}}
\newcommand{\ind}{\mathcal{J}}

\newcommand{\La}{\Lambda}
\newcommand{\Nr}{\bn_0^r}

\newcommand{\Toep}{\mathcal{T}_{\Lambda}}
\newcommand{\Har}{\mathcal{A}_{\Lambda}}
\newcommand{\wHar}{\mathcal{H}_{\Lambda}}
\newcommand{\Cun}{\mathcal{O}_{\Lambda}}
\newcommand{\Gauge}{\mathcal{F}_{\Lambda}}

\newcommand{\Lin}{\textrm{Lin}}

\newcommand{\LFock}{l^2(\Lambda)}

\newcommand{\oner}{\{1,\ldots,r\}}

\newenvironment{rlist}
{

\begin{enumerate}}
{\end{enumerate}}

\newcommand{\la}{\lambda}
\newcommand{\ra}{\rangle}
\newcommand{\rra}{\right\rangle}
\newcommand{\lla}{\left\langle}

\newcommand{\ot}{\otimes}

\newcommand{\ol}{\overline}

\newcommand{\wt}{\widetilde}

\numberwithin{equation}{section}


\begin{document}
\author{Adam Skalski}
\author{Joachim Zacharias}
\footnote{\emph{Permanent address of the first named author:}
Department of Mathematics, University of \L\'{o}d\'{z}, ul.
Banacha 22, 90-238 \L\'{o}d\'{z}, Poland.

Research supported by the EPSRC grant no. RIS 24893}
%\footnote{\emph{Permanent address of the first named author:}
%Department of Mathematics, University of \L\'{o}d\'{z}, ul.
%Banacha 22, 90-238 \L\'{o}d\'{z}, Poland.}

\address{School of Mathematical Sciences,  University of Nottingham,
Nottingham, NG7 2RD} \email{adam.skalski@maths.nottingham.ac.uk} \email{joachim.zacharias@nottingham.ac.uk  }

\title{\bf Poisson transform for higher-rank graph algebras and its applications}
\begin{abstract}
Higher-rank graph generalisations of the Popescu-Poisson transform are constructed, allowing us to develop a dilation theory for
higher rank operator tuples. These dilations are joint dilations of the families of operators satisfying relations encoded by the
graph structure which we call $\Lambda$-contractions or $\Lambda$-isometries. Besides commutant lifting results and
characterisations of pure states on higher rank graph algebras several applications to the structure theory of non-selfadjoint
graph operator algebras are presented generalising recent results in special cases.
\end{abstract}

\maketitle
The starting point of classical dilation theory was the construction of a dilation of a
contractive operator on a Hilbert space to a unitary by B.\,Sz.-Nagy in 1953 (\cite{Nagy}). Under a certain minimality
condition this dilation is unique up to unitary equivalence.
There exists also a minimal isometric
dilation whose adjoint leaves the original Hilbert space invariant.
%\cite{Agler}
%(we suggest the term coextension).
In particular, this provides a quick proof of von Neumann's celebrated inequality for Hilbert space contractions.

%Another basic dilation theorem is of course Stinespring's Theorem showing that a completely positive
%maps from a $C^*$-algebra to $B(\Hil)$ allows dilations to $C^*$-representations again
%unique under certain minimality conditions.

There is a connection between Sz.-Nagy's and Stinespring's Theorems: a contraction defines
a completely positive map and using this it turns out that the existence of minimal isometric dilations
of contractions can be derived from Stinespring's theorem. This idea links dilation theory of contractions to
representations of the Toeplitz algebra and plays a key role in the more recent developments of dilation
theory for tuples of operators.

In 1989 G.\,Popescu proved  in \cite{Popescu} a result  analogous to Sz.-Nagy's for a \emph{row-contraction}, that is a (finite
or infinite) sequence of contractions $(T_i)_{i=1}^n$ on a Hilbert space $\Hil$ such that $\sum_{i=1}^n T_i T_i^* \leq I_{\Hil}$.
Building on the earlier work of A.E.\,Frazho and J.\,Bunce he showed that each such object can be dilated to a
\emph{row-isometry}, which in turn provides a representation of the Cuntz algebra $\mathcal{O}_n$. This may be regarded as a
dilation of the completely positive map $X \mapsto \sum_{i=1}^n T_iXT_i^*$. In a series of further papers G.\,Popescu unveiled
the whole array of connections between row-contractions, their dilations, operator algebras related to free semigroups on
$n$-generators, completely positive maps on $B(\Hil)$ and their dilations to endomorphisms. This led him to develop a theory
which can be justifiably called noncommutative (free) complex analysis (see for example \cite{Pcomp} and references therein),
with concrete operators on a full Fock space playing the role of certain classes of analytic functions; in particular, it allowed
to establish a von Neumann inequality for row contractions.
His most important tool is a noncommutative Poisson transform %which
%links the dilation of tuples with Stinespring's Theorem and representations of Cuntz algebras
(introduced in \cite{PPois} and later used in \cite{pureod}).

This $C^*$-algebraic formulation of dilation theory for row contractions has subsequently also been used to construct dilations
of tuples preserving symmetry conditions. In \cite{subalgIII} W.\,Arveson constructed a universal commuting row contraction
(given by shifts on the symmetric Fock space) and a corresponding Toeplitz $C^*$-algebra. He studied dilations of commuting row
contractions via dilating the corresponding completely positive maps to representations of this Toeplitz algebra. Similar
theories have been developed for $q$-commuting tuples by S.\,Dey (\cite{Santanu}) and for tuples verifying conditions analogous
to these satisfied by the generators of Cuntz-Krieger algebras by B.V.R.\,Bhat, S.\,Dey and the second named author
(\cite{RajJSan}).

Each algebra occuring in these dilation theories is generated by a family of operators forming a single row contraction.
A very general class of algebras generated in this way is the class of graph $C^*$-algebras which has been studied intensively
in recent years.
Initially inspired by the Cuntz-Krieger algebras, it was
soon shown to provide a rich source of new examples of $C^*$-algebras and also allowed to view and
analyse certain known algebras from a new perspective.
An even more general class is the one of higher rank graph algebras introduced by A.\,Kumjian and D.\,Pask in \cite{kupa}. These were inspired by the
higher rank Cuntz-Krieger algebras of G.\,Robertson and T.\,Steger (\cite{Steg}).
Graph algebras are usually
defined as certain universal objects or algebras related to groupoids.
%For the purposes of this paper we prefer to work with the concrete algebras of operators
% acting on .
Note however that exactly as the Cuntz algebra $\mathcal{O}_n$ is a quotient of the concrete
Toeplitz-type algebra acting on the full Fock space of
$\bc^n$, the same remains true for its graph generalisations; the graph Cuntz-type algebra $\Cun$ arises as a
quotient of a certain concrete algebra of operators $\Toep$, acting on the $L^2$-space of the higher rank graph $\La$.
%a only that the ideal is more complicated.
In this picture $\mathcal{O}_n$ arises from a graph
which has one vertex and $n$ edges. %A precise definition and comments on the connections with the universal constructions are given in Section 1.
For an exhaustive outline of the theory of the graph algebras we refer to the book \cite{book} and references therein.
Note that a rank $r$-graph algebra ($r \in \bn$) can be thought of as an algebra generated by $r$  row contractions
with specific commutation relations encoded by the higher rank graph. The relevant generators may thus be regarded as
higher rank tuples or higher rank row contractions.

%Here we would like to point out
%that exactly as the Cuntz algebra
%$\mathcal{O}_n$ is equal to a quotient of the concrete Toeplitz-type algebras acting on the full Fock space of
%$\bc^n$, the same remains true for its graph generalisation, with the Fock space replaced with the $L^2$-space over
%paths in the graph under consideration. In this picture $\mathcal{O}_n$ arises from a graph which has one vertex
%and $n$ edges. A precise definition and comments on the connections with the universal constructions are given in
%Section 1. We work with the higher-rank graphs, which were introduced in \cite{kupa} and generalise standard
%(rank-1) graphs.

The aim of this paper is to develop a dilation theory in the framework of higher rank tuples. To stress the dependence on $\La$
we call the basic objects of our dilation theory \emph{$\La$-contractions}. $\La$-contractions are families of contractions on a
Hilbert space $\Hil$ indexed by paths of $\La$ encoding the relevant graph structure. The set of vertices of $\La$ is
understood to induce the decomposition of $\Hil$ into orthogonal subspaces, and an operator corresponding to a path $\lambda$ is
assumed to act nontrivially only between the subspaces corresponding to the source and range of $\lambda$. The additional
constraints correspond to the row-contraction-type condition. If the rank of the graph is greater than one, we may think of edges
as having various colours. Then the conditions defining $\La$-contractions encode also certain commutation relations between the
operators corresponding to edges of varying colours. In the rank-1 case certain dilation results for objects of that type were
established by M.T.\,Jury and D.W.\,Kribs in \cite{Jury}. Their starting point was however different: putting it in our language
they consider a given row contraction and then look at the class of potential graphs $\La$ for which the row contraction in
question may be viewed as a $\La$-contraction (the procedure involves incorporating for each $\La$ an appropriate family of
vertex projections). This led them to introduce a certain partial ordering on the class of resulting dilations. It is not clear
how to extend their approach to higher-rank cases.

The basic tool for the analysis of $\La$-contractions will be a generalisation of the Popescu-Poisson transform introduced in
\cite{PPois}. Already in that paper the transform is constructed for certain higher-rank objects but restricted to tensor
products of standard Cuntz (respectively Cuntz-Toeplitz) algebras. As a $\La$-contraction (of rank $r$) allows to define $r$
commuting completely positive contractions, thus defining a canonical semigroup of completely positive maps ($\sigma_{\Vcont}$),
our dilation problem amounts to finding a particular dilation of this semigroup. It is known from the classical theory that one
cannot expect joint dilations of three (or more) commuting contractions to hold without any additional assumptions.  This
explains why it was necessary in \cite{PPois} to restrict attention to the families of operators satisfying the condition $(P)$.
We introduce an analogous constraint for $\La$-contractions and call it the \emph{Popescu condition} (see Definition
\ref{Popcond}). It is shown that every $\La$-contraction satisfying the Popescu condition admits a dilation to a
\emph{Toeplitz-Cuntz-Krieger} family, which is essentially a family of contractions arising as a representation of the
Toeplitz-type algebra $\Toep$. Similarly each $\La$-isometry  admits a dilation to a \emph{Cuntz-Pimsner} family, a family of
contractions arising as a representation of the Cuntz-type algebra $\Cun$. These dilations
%especially in the higher-rank case,
may be thought of as joint dilations of (higher rank) tuples
satisfying certain commutation relations in such a way that the commutation relations are preserved.
%uch constructions in particular cases we refer to \cite{BBD} and \cite{RajJSan}.
Graph versions of the Popescu-Poisson transform also allow to establish two types of commutant lifting theorems, one related to the
commutant of a given $\La$-contraction $\Vcont$ and another to the fixed point space of the semigroup of
completely positive maps associated with $\Vcont$. Following the methods of G.\,Popescu, K.\,Davidson
and D.\,Pitts used for the noncommutative $H^{\infty}$ algebras associated to free semigroups (\cite{Pderiv} and
\cite{DavPit}), we can further compute the character spaces of the  Hardy-type Banach algebras $\Har$ and $\wHar$
and show that $\Har$ is never amenable (when nontrivial). Finally the GNS construction relative to a state on
$\Cun$ is given an alternative description in terms of `cyclic' $\La$-isometries, and, as was done in
\cite{pureod} for the standard Cuntz algebras, purity of the state is characterised in terms of the ergodicity of the
associated semigroup of completely positive maps.

The detailed plan of the paper is as follows. In Section 1 we recall the definition of higher-rank graphs,
introduce their associated concrete operator algebras (and comment on relations with the universal ones),
define $\La$-contractions, $\La$-isometries and the Popescu condition. Section 2 contains the definition
and basic properties of the Toeplitz-Poisson and Cuntz-Poisson transforms including a von Neumann inequality for
$\La$-contractions. In Section 3 we show that the minimal Stinespring construction for the Toeplitz-Poisson transform
associated with a given $\La$-contraction $\Vcont$ yields a minimal dilation of $\Vcont$ to a Toeplitz-Cuntz-Krieger family and that if $\Vcont$ is a $\La$-isometry then the dilation automatically forms a Cuntz-Pimsner family.
Our two different commutant lifting theorems are also established in this section.
%for the commutant of $\Vcont$ and for the fixed point space of the canonical semigroup $\sigma_{\Vcont}$.
Finally Section 4 contains a number of applications of our Poisson-type transforms.
We use them to describe the character spaces of the algebras $\Har$ and $\wHar$, to prove nonamenability of $\Har$
for all nontrivial graphs $\La$ and to characterise pure states on $\Cun$ in terms of the related $\La$-isometry.

%Note that in the introduction we used the symbol $\La$ to denote a graph, whereas further on, following the established
%recently terminology in the operator algebra theory,
%it will denote a category of paths associated with a fixed (higher-rank) graph. We hope it will not cause any confusion.

\section{Notations and basic definitions}

This section begins with recalling the definition of higher-rank graphs,
Toeplitz-Cuntz-Krieger and Cuntz-Pimsner families and related concrete graph operator algebras.
We proceed to introduce the crucial concept of $\La$-contractions and $\La$-isometries and discuss the
Popescu condition.

\subsection*{Higher-rank graphs}

Higher rank graphs were introduced in \cite{kupa} as a tool for constructing $C^*$-algebras generalising higher
rank Cuntz-Krieger algebras defined in \cite{Steg}.  In this paper we mostly follow the notation in \cite{book} some of
which we recall briefly.

Let $\bn_0= \bn \cup \{0\}$. For $r\in \bn$ the canonical `basis' in $\Nr$ will be denoted by $(e_1, \ldots, e_r)$, and we write $e=\sum_{i=1}^r
e_i$. The componentwise maximum of $n,m \in \Nr$ is denoted by $n \vee m$ and we write
 $|n|= n_1 + \cdots + n_r$. Throughout the paper $\La$ will denote a rank-$r$ graph, that is a small category
with set of objects $\Lambda^0$ and shape functor $\sigma : \Lambda \to {\bf N}^r$
(where ${\bf N}^r$ is viewed as the  category with one object and
morphisms $\Nr $) satisfying the \emph{factorisation property} defined in \cite{kupa}. If $n \in \Nr$ the set of morphisms
in $\La$ of shape $n$ is denoted by $\La^n$. For $\lambda \in \Lambda$ we write
$|\lambda |=|\sigma(\lambda)|$. The morphisms in $\Lambda$
may be thought of as paths in a `multi-coloured' graph with vertices indexed by the set $\La^0$.
The range and source maps are respectively denoted by $r:\La \to \La^0$ and $s:\La \to \La^0$
(called origin and terminus in \cite{Steg}). The factorisation property says that if $m,n \in \Nr$ then
every morphism $\lambda \in \Lambda^{m+n}$ is a unique product $\lambda=\mu \nu$ of a  $\mu \in \Lambda^m$ and
$\nu \in \Lambda^n$, where $s(\mu)=r(\nu)$. To simplify the language we will assume that $\La$ is countable.
In fact all the results remain valid if this assumption is dropped.

A rank-$r$ graph $\La$ is called \emph{finitely aligned} if for each $\la,\mu \in \La$ the set of minimal common extensions of
$\la$ and $\mu$, that is $MCE(\la, \mu):=\{ \nu \in \La : \exists_{\alpha ,\beta \in \Lambda} \; \nu = \la \alpha = \mu \beta, \:
\sigma(\la \alpha) = \sigma(\la) \vee \sigma(\mu)\}$, is finite. It is called \emph{row-finite} if for each $a \in \La^0$ and $n
\in \Nr$ the set $\La^n_a:=\{ \la \in \La: r(\la) = a, \sigma(\la) =n\}$ of paths of shape $n$ ending at $a$ is finite. When $a,b
\in \La^0$, $n \in \Nr$, we also write $\La^n_{b,a}:=\{ \la \in \La: r(\la) = a, s (\la)=b, \sigma(\la) =n\}$. Finally we call
graph $\La$ \emph{finite} if $\La^n$ is finite for any $n \in \Nr$ (note that it does not mean that $\La$ is finite as a set). It
is immediate that finite graphs are row-finite, and row-finite graphs are finitely aligned. $\La$ is said to have \emph{no
sources} if for each $n \in \Nr$ and $a \in \La^0$ there exists $\la \in \La^n$ such that $r(\la)=a$, equivalently $\La_a^n \neq
\emptyset$ (there are arbitrarily long paths ending at $a$). A weaker condition is the following: $\La$ is said to be
\emph{cofinal} if for every $a \in \La^0$ there is $\la \in \La$ with $\sigma (\la) \neq 0$ and $r(\la)=a$.

The following definition was introduced in \cite{toep}. Note that we use a different convention to the one
used in that paper when it comes to labelling target and initial projections of the partial
isometries in question (to maintain compatibility with \cite{book} source and range are exchanged).
\begin{deft}\label{TCK}
Suppose that $\La$ is finitely aligned.
A family of partial isometries $\{x_{\la}:\la \in \La\}$ in a $C^*$-algebra $B$
is called a Toeplitz-Cuntz-Krieger $\La$-family if the following are satisfied:
\begin{rlist}
\item $\{x_a: a \in \La^0\}$ is a family of mutually orthogonal projections;
\item $x_{\la} x_{\mu} = x_{\la \mu}$ if $\la, \mu \in \La$, $s(\la) = r(\mu)$;
\item $x_{\la}^* x_{\la} = x_{s(\la)}$ if $\la \in \La$;
\item if $n \in \Nr\setminus\{0\}$, $a \in \La^0$  and $F\subset \La^n_a$ is finite then
$x_a \geq \sum_{\la \in F} x_{\la} x_{\la}^*$;
\item $x_{\mu}^* x _{\nu} = \sum_{\mu \alpha = \nu \beta \in MCE(\mu, \nu)} x_{\alpha} x_{\beta}^*$
for all $\mu, \nu \in \La$.
\end{rlist}
% The sum in (iv) is to be understood in the strong operator topology.
If $\la$ is row finite and instead of (iv)  a stronger condition
\[x_a = \sum_{\la \in \La^n_a} x_{\la} x_{\la}^*\]
is satisfied for all $n \in \Nr\setminus\{0\}$ and $a \in \La^0$ then $\{x_{\la}:\la \in \La\}$ is called a Cuntz-Pimsner
$\La$-family.
\end{deft}

\subsection*{Concrete operator algebras associated with $\La$}

Consider the creation operators on a Hilbert space $l^2(\Lambda)$
(where $\La$ is equipped with the counting measure) defined by
\[ L_{\lambda} \delta_{\mu} = \left\{\begin{array}{ccc} \delta_{\lambda \mu}& \textrm{if} & s(\lambda) =  r(\mu),\\
0& \textrm{if} & s(\lambda)\neq r(\mu),\end{array}, \;\;\; \lambda, \mu \in \Lambda. \right.
\]
 Further we follow the convention that $\delta_{\lambda \mu}:=0$ and $L_{\lambda \mu }:=0$ if
$s(\lambda)\neq r(\mu)$. It is easily seen that all $L_{\lambda}$ are partial isometries. If $\La$ is finitely
aligned they form a Toeplitz-Cuntz-Krieger family (in $B(l^2(\La)))$.
The space $l^2(\Lambda)$ may be viewed as a Fock space associated with the graph $\La$. In particular
if $r=1$ and the graph $\La$ has only one vertex with $k$ edges, $\LFock$ is naturally isometrically isomorphic
to a free Fock space over $\bc^k$.

We will consider the following concrete operator algebras contained in $B(\LFock)$:
\[ \Toep = C^*(\{L_{\lambda}:\lambda \in \La\}),\]
\[ \Har = \mathcal{ALG}(\{I, L_{\lambda}:\lambda \in \La\}),\]
\[\wHar = \mathcal{ALG}_w(\{I, L_{\lambda}:\lambda \in \La\}),\]
where $\mathcal{ALG}(S)$ (respectively, $\mathcal{ALG}_w(S)$) denotes the smallest norm (respectively ultraweakly)
closed subalgebra generated by the set $S$.
The algebras above can be respectively interpreted as noncommutative higher-rank graph versions of the Toeplitz
$C^*$-algebra, disc algebra, and $H^{\infty}$.
Note that $\Toep$ is unital if and only if the set $\La^0$ is finite.
If this is not the case we will occasionally consider the natural unitisation $\Toep^1:=C^*(I, \Toep)$.
It is important in the sequel to note that $\Toep= \overline{\Lin}
\{L_{\mu} L_{\nu}^*:\mu, \nu \in \La\}$ when $\La$ is finitely aligned.
This is an immediate consequence of  equality (iii) in Definition \ref{TCK}.

Denote by $P_j$ ($j \in \oner$) the projection onto the span of those $\delta_{\lambda}$ for which
$\sigma(\lambda)_j=0$. Note that if $\La$ is finite then
\[ P_j = I - \sum_{\la \in \La^{e_j}} L_{\la} L_{\la}^* =  \sum_{a \in \La^0}L_a  - \sum_{\la \in \La^{e_j}} L_{\la} L_{\la}^* ,\]
so that each $P_j$ belongs to $\Toep$. In this case we denote by $\mathcal{J}$ the ideal of $\Toep$
generated by all $P_j$ and define the Cuntz-type algebra
\[ \Cun = \Toep /\mathcal{J}.\]
The canonical generators in $\Cun$, obtained by quotienting from the creation operators $L_{\la}$, will be denoted
by $s_{\la}$. These form a Cuntz-Pimsner $\La$-family.

It is natural to ask whether the Toeplitz- and Cuntz-type algebras introduced above coincide with `universal'
graph algebras considered in \cite{kupa}, \cite{toep} and \cite{align}. Recall that the latter are defined as
universal objects respectively for Toeplitz-Cuntz-Krieger and Cuntz-Pimsner  families.
When $\La$ is finitely aligned, Corollary
7.7 of \cite{toep} implies that $\Toep$ coincides with the corresponding universal object. For the Cuntz-type algebra
the situation is a little more complicated, but if $\La$ is finite and has no sources one can exploit in the usual way
the gauge uniqueness theorem. Define first for each $z \in \bt^r$ a unitary $U_z \in B(l^2(\La))$ by the linear continuous extension
of the formula
\[U_z (\delta_\la) = z^{\sigma(\la)} \delta_{\la}, \;\;\;\;\la \in \La,\]
where for all $z \in \bt^r, n \in \Nr$ we write $z^n:= z_1^{n_1} \cdots z_r^{n_r}$. The formula
\[\alpha_z(x) = U_z x U_z^*, \;\;\; x \in \Toep\]
is easily seen to define an automorphism of $\Toep$. As we have $\alpha_z(L_{\la})= z^{\sigma(\la)} L_{\la}$
for each $\la \in \La$, it is also easy to see that $\alpha$ yields
a (pointwise norm) continuous action of $\bt^r$ on $\Toep$. For $\La$ finite each $\alpha_z$ leaves the ideal $\mathcal J$
invariant, so that the action of $\alpha$ descends canonically to $\Cun$. Now it remains to observe that if $\La$
has no sources then each of the generators $s_a$ ($a \in \La^0$) is non-zero. Indeed, it is easy to see that
whenever $\mu, \nu, \alpha, \beta \in \La$, and $j \in \oner$ then for all $\gamma \in \La$ such that
$\sigma(\gamma)_j > \sigma(\beta)_j - \sigma(\alpha)_j$
\[L_{\mu} L_{\nu}^* P_j L_{\alpha} L_{\beta}^* \delta_{\gamma} = 0.\]
On the other hand for any $n \in \Nr$ we can find $\gamma \in \La^n_a$ such that $L_a \delta_{\gamma}  = \delta_{\gamma}$. As the
set of finite linear combinations of operators of the form $L_{\mu} L_{\nu}^* P_j L_{\alpha} L_{\beta}^*$ is norm dense in
$\mathcal{J}$, we see that $d(L_a, \mathcal{J})=1$.
 Theorem 3.4 of \cite{kupa}
implies that $\Cun$ is isomorphic to the universal object $C^*(\La)$.

It is clear from the discussion above that $\mathcal{J}$ may be viewed as the ideal of all operators of
`finite length in at least one direction' in $\Toep$. If $\La^0$ is infinite, $\Toep$ may not contain
any operators of that type. Providing Fock space models for $\Cun$ becomes more complicated in this situation.

\subsection*{$\La$-contractions and $\La$-isometries}

In Section 3 we will be interested in joint dilations of families of operators to Toeplitz-Cuntz-Krieger
or Cuntz-Pimsner families. The definitions below encapsulate the conditions we would impose on the families
we expect to dilate.

\begin{deft}        \label{Lcont}
Let $\Hil$ be a Hilbert space. A family $\Vcont =\{V_{\lambda}:\la \in \La\}$ of operators in $B(\Hil)$ is called
a $\La$-contraction if the following conditions are satisfied:
\begin{rlist}
\item $\forall_{\la,\mu \in \Lambda,\, s(\la)\neq r(\mu)} \; V_{\la} V_{\mu} =0$;
\item $\forall_{\la,\mu \in \Lambda,\, s(\la) = r(\mu)} \; V_{\la} V_{\mu} =V_{\la \mu}$;
\item $\forall_{n \in \Nr} \; \sum_{\la \in \La^n} V_{\la} V_{\la}^* \leq I$;
\item each $V_a$ ($a \in \La^0$) is an orthogonal projection and $\sum_{a\in \La^0} V_a =I$
\end{rlist}
All infinite sums of Hilbert space operators here and in what follows are understood in the strong operator topology.
\end{deft}

The definition above is not optimal -  condition (iv) is in a sense redundant. If $\Vcont =\{V_{\lambda}:\la \in
\La\}$ is a family of operators satisfying (i)-(iii) above then conditions  (ii) and (iii) imply that each $V_a$
for $a \in \La^0$ is a contractive idempotent, hence a projection. Further (i) shows that $V_a V_b = 0$ if $b \in
\La^0$ and $a\neq b$. Denoting by $p$ the sum $\sum_{a \in \La^0} V_a$ we see that $V_{\la} =
pV_{\la}p$ (by (i) and (ii)). Therefore even if the second part of the condition (iv) is not satisfied at
the outset, we can to all aims and purposes analyse the obvious $\La$-contraction on $p\Hil$. Following the
convention introduced earlier we define $V_{\lambda \mu}:=0$ if $s(\la) \neq r(\mu)$.

\begin{deft}         \label{Lisom}
A $\La$-contraction is called a $\La$-isometry if the following holds:
\begin{rlist}\setcounter{enumi}{4}
\item $\forall_{n \in \Nr} \; \sum_{\la \in \La^n} V_{\la} V_{\la}^* = I$.
\end{rlist}
\end{deft}

Following on the theme of non-optimality of the definitions above, in conditions (iii) and (v) above it is enough
to assume that the inequalities (respectively, equalities) hold only for $n$ of the form $e_j$, $j \in \oner$.

We point out that the existence of nontrivial $\La$-isometries imposes conditions on the graph $\La$.
Indeed, $\sum_{\la \in \La^n} V_{\la} V_{\la}^* = I$
for $n \in \Nr$ implies $V_a= V_a \sum_{\la \in \La^n} V_{\la} V_{\la}^* =   \sum_{\la \in \La_a^n} V_{\la} V_{\la}^*$
so that if only $V_a$ is nonzero, $\La_a^n \neq \emptyset$. This implies that if $\Hil$ is nontrivial, then
$\La$ has to be cofinal, and if each of the vertex projections $V_a$ is non-zero then $\La$ has no sources.
Conversely, if $\La$ has no sources then nontrivial $\La$-isometries exist
(just consider the canonical Cuntz-Pimsner $\La$-family in $B(l^2(\La))$.)

One may think of $\Lambda$-contractions as generated by $r$ sets of row-contractions, each corresponding
to a `coordinate' of $\La$, with certain commutation relations imposed on them. In particular when the graph $\La$
is finite, a $\Lambda$-contraction is generated by a finite set of operators satisfying a combination of
inequalities and commutation relations. This is related to the way of seeing higher-rank graphs as product systems
of usual graphs over $\oner$, as analysed by I.\,Raeburn and A.\,Sims in \cite{toep}.

\begin{deft} \label{Popcond}
Let $\Vcont$ be a $\La$-contraction and define for $s \in (0,1)$ the defect operator
\begin{equation}\label{Delta} \Pop =
\sum_{\mu \in \La, \, \sigma(\mu) \leq e} (-s^2)^{|\mu|} V_{\mu} V_{\mu}^*.
\end{equation}
The family $\Vcont$ is said to satisfy the Popescu condition (or condition `P') if there exists $\rho\in (0,1)$
such that for all $s \in (\rho,1)$ the operator $\Pop$ is positive.
\end{deft}

Each $\La$-isometry $\Vcont$ satisfies the Popescu condition, as then for each $s\in (0,1)$ one can easily compute $\Pop=
\sum_{\mu \in \La, \, \sigma(\mu) \leq e} (-s^2)^{|\mu|} I = (1-s^2)^r I$. Another case when the Popescu condition is
automatically satisfied is the one of \emph{doubly commuting} $\La$-contractions. Here a $\La$-contraction $\Vcont$ is doubly
commuting if and only if for all $i,j\in \oner$, $\la \in \La^{e_i}$, $\mu \in \La^{e_j}$, there is
\[ V_{\la}^* V_{\mu} = \sum_{\alpha \in \La^{e_i}, \beta \in \La^{e_j}, \mu \beta = \la \alpha} V_{\beta}
V_{\alpha}^*.\] Indeed, for a doubly commuting $\La$-contraction $\Vcont$
\[ \Pop = \prod_{j=1}^r (I - s^2 \sum_{\mu \in \La^{e_j}} V_{\mu}
V_{\mu}^*),\] and the right hand side of the above expression is positive as the product of commuting positive
operators for all $s \in (0,1)$.

\section{The Poisson transform for $\Lambda$-contractions}
In this section we define and introduce basic properties of the Poisson transform associated to a higher-rank graph. The
techniques of the proofs are similar to the ones in \cite{PPois}, but the graph formalism not only yields a far-reaching
generalisation of the results of that paper but also allows for  more concise, and (in our opinion) more transparent proofs. The
theorems established here provide methods to be applied in the following sections.

\begin{lem}
Let $\Vcont$ be a $\La$-contraction (on $\Hil$) and recall the operator $\Pop$ defined by \eqref{Delta}. Then for
any $s \in (0,1)$ the following holds
\begin{equation} \label{absorb}
\sum_{\la \in \La} s^{2|\la|} V_{\la}  \Pop V_{\la}^*  = I_{\Hil}
\end{equation}
In particular if %$\Vcont$ satisfies the Popescu condition and  $s \in (0,1)$ is such that
the operator $\Pop$  is
positive then the map $\Wop \in B(\Hil;\LFock \ot \Hil)$ defined by
\begin{equation} \label{Welta}
\Wop (\xi) = \sum_{\la \in \La} \delta_{\la} \ot s^{|\la|} \Pop^{\frac{1}{2}} V_{\la}^* \xi
\end{equation}
($\xi \in \Hil$) is an isometry.
\end{lem}

\begin{proof}
Compute:
\begin{align*} \sum_{\la \in \La} s^{2|\la|} V_{\la}
\Pop V_{\la}^*&
  = \sum_{\la \in \La} s^{2|\la|} V_{\la} \left( \sum_{\mu \in \La, \, \sigma(\mu) \leq e}
  (-s^2)^{|\mu|} V_{\mu} V_{\mu}^* \right)   V_{\la}^*  \\
  =& \sum_{\la \in \La} \sum_{\mu \in \La, \, \sigma(\mu) \leq e} s^{2(|\la|+ |\mu|)}
   (-1)^{|\mu|} V_{\la \mu} V_{\la\mu}^*
  = \sum_{\gamma \in \La}  s^{2|\gamma|} C_{\gamma} V_{\gamma} V_{\gamma}^* ,
\end{align*}
where $C_{\gamma} \in \br, (\gamma \in \La)$ are certain constants. To determine $C_{\gamma}$ we have to consider all
factorisations $\gamma = \alpha \mu$ with $\sigma (\mu ) \leq e$. If $\sigma(\gamma) = 0$, then the only way $\gamma$ can be
written as the composition of two elements is $\gamma= \gamma \gamma$, so $C_{\gamma} = (-1)^0 = 1$.

Suppose now that for some $k \in \oner$ and some non-empty set $\mathcal{I} =\{i_1,\ldots,i_k\} \subset
\oner$ of indices $\sigma(\gamma)_i>0$ for $i \in \mathcal{I}$ and $\sigma(\gamma)_j =0$ for $j \in \oner \setminus
\mathcal{I}$. Let $m \in \Nr$, $m \leq e$. Then $\gamma$ can be written (in a unique way) as $\alpha \mu$ for some
$\alpha \in \La$ and $\mu \in \La^m$
 if and only if $m \leq e_{i_1} + \cdots +e_{i_k}$. This implies
 that
\[C_{\gamma} = \sum_{i_1,\ldots,i_k=0}^1 (-1)^{i_1 + \cdots + i_k}=
0.\]

Combination of the above shows that
\[
 \sum_{\la \in \La} s^{2|\la|}
V_{\la} \Pop V_{\la}^*
  = \sum_{\gamma \in \La^0}  V_{\gamma} V_{\gamma}^*  = I_{\Hil}.
\]
%The vector $\xi$ obviously plays no role in the above computation - we write it explicitly to stress that all infinite
% sums of operators, as usual, are to be understood in the strong operator topology.

The second part of the lemma follows from the first.
\end{proof}

With the above in hand we are ready to establish the main theorem of this section.
% Note that for the uniqueness there we need to assume that $\La$ is finitely aligned.

\begin{tw} \label{tPoisson}
Let $\La$ be finitely aligned and let $\Vcont$ be a $\La$-contraction (on $\Hil$) satisfying the Popescu condition. Then there
exists a unique continuous linear map $R_{\Vcont}:\Toep \to B(\Hil)$ satisfying
\[R_{\Vcont}(L_{\la} L_{\mu}^*) = V_{\la} V_{\mu}^*, \;\;\;\; \la, \mu \in \La.\]
The map $R_{\Vcont}$ will be called the $\La$-Poisson transform (associated
with $\Vcont$). It is completely positive and contractive, unital if $\Toep$ is unital. If $\Toep$ is not unital,
$R_{\Vcont}$ has a unique unital extension to $\Toep^1$, which is also completely positive. In any case
$R_{\Vcont}|_{\Har}$ is multiplicative.
\end{tw}

\begin{proof}
Let $\rho\in (0,1)$ be such that for all $s \in (\rho,1)$ the operator $\Pop$ is positive. Define for each such $s$ a
map $\sPois:\Toep \to B(\Hil)$ by the formula
\[\sPois (x) = \Wop^* (x \ot I_{\Hil}) \Wop,\]
where $\Wop$ is the isometry defined in the last lemma (note that the above formula obviously makes sense for any $x \in
B(\LFock)$). It is clear that $\sPois$ is completely positive and contractive. Moreover for any $\mu,\nu \in \La$
\begin{equation} \label{rPois}
\sPois (L_{\mu} L_{\nu}^*) = s^{|\mu| + |\nu|}V_{\mu} V_{\nu}^*.
\end{equation}
Indeed, let $\mu,\nu \in \La$ and $\xi, \eta \in \Hil$. Then, using the convention $V_{\alpha \beta} :=0$
if $s(\alpha) \neq r(\beta)$, we have
\begin{align*}
\lla \eta, \sPois (L_{\mu} L_{\nu}^*) \xi \rra &= \lla \Wop
\eta,  \left( L_{\mu} L_{\nu}^* \ot I_{\Hil}  \right)\Wop \xi \rra \\
&= \lla \sum_{\la \in \La} \delta_{\la} \ot s^{|\la|} \Pop^{\frac{1}{2}} V_{\la}^* \eta, \left(L_{\mu} L_{\nu}^*
\ot I_{\Hil} \right) \sum_{\gamma \in \La} \delta_{\gamma} \ot s^{|\gamma|}
\Pop^{\frac{1}{2}} V_{\gamma}^* \xi \rra \\
&=
 \lla \sum_{\la \in \La} \delta_{\la} \ot
s^{|\la|} \Pop^{\frac{1}{2}} V_{\la}^* \eta,  \sum_{\alpha \in \La} \delta_{\mu\alpha} \ot s^{|\alpha|+|\nu|}
\Pop^{\frac{1}{2}}
V_{\nu\alpha}^* \xi \rra  \\
&= \lla \sum_{\beta \in \La} \delta_{\mu\beta} \ot s^{|\mu|+|\beta|} \Pop^{\frac{1}{2}} V_{\mu \beta}^* \eta,
\sum_{\alpha \in \La} \delta_{\mu\alpha} \ot s^{|\alpha|+|\nu|} \Pop^{\frac{1}{2}}
V_{\nu\alpha}^* \xi \rra \\
&=
  s^{|\mu|+|\nu|} \sum_{\alpha \in \La}  s^{2 |\alpha|}
 \lla \Pop^{\frac{1}{2}} V_{\mu
\alpha}^* \eta,  \Pop^{\frac{1}{2}} V_{\nu\alpha}^* \xi \rra \\&=
 s^{|\mu|+|\nu|} \lla V_{\mu}^* \eta,
  \left( \sum_{\alpha \in \La}  s^{2 |\alpha|}
 V_{\alpha}^* \Pop V_{\alpha} \right) V_{\nu}^* \xi \rra.
\end{align*}
Now \eqref{absorb} implies that the above is actually equal to $ s^{|\mu|+|\nu|} \langle V_{\mu}^* \eta,
   V_{\nu}^* \xi \ra,$
and the equality \eqref{rPois} is established. As the set $\{L_{\mu} L_{\nu}^*: \mu, \nu \in \La\}$ is total in
$\Toep$, it follows by standard arguments that for each $x \in \Toep$ the limit $\lim_{s \to 1^-} \sPois (x)$
exists (in the norm topology), and moreover the map $R_{\Vcont}:\Toep \to B(\Hil)$ defined by
\[ R_{\Vcont}(x) = \lim_{s
\to 1^-} \sPois (x), \;\;\; x \in \Toep\] satisfies all the requirements of the theorem. Uniqueness is again a
consequence of the totality of the set $\{L_{\mu} L_{\nu}^*: \mu, \nu \in \La\}$ in $\Toep$.
All the statements concerning unitality and multiplicativity are now easy to prove.
\end{proof}

\begin{rem} \label{notal}
Note that if we are only interested in the operator algebra $\Har$, to obtain the existence of the completely contractive
map $R_{\Vcont}: \Har \to B(\Hil)$ with all the properties above we do not need to assume that $\La$ is finitely aligned. The latter condition is only necessary to prove uniqueness in the $C^*$-case.
\end{rem}

Suppose again that $\La$ is finitely aligned. For each operator polynomial $p({\mathcal L}) \in B(\ell^2(\Lambda))$ in noncommuting variables $L_{\mu}$ and $L_{\nu}^*$ there exist finitely many complex coefficients
$\{\alpha_{\kappa, \nu} \in \bc: \kappa, \nu \in \La\}$ such that
\begin{equation}
p({\mathcal L})= \sum_{\kappa, \nu \in \La} \alpha_{\kappa, \nu}
  L_{\kappa} L_{\nu}^*.
\label{poly}
\end{equation}
%(here ${\mathcal L}$ stands for the $\Lambda$-contraction ${\mathcal L}=\{L_{\la}: \la \in \La\} $).
Following
\cite{PPois} for any $\La$-contraction on a Hilbert space $\Hil$ we define the operator $p (\Vcont)\in B(\Hil)$ by
\[ p (\Vcont) = \sum_{\kappa, \nu \in \La} \alpha_{\kappa, \nu}
  V_{\kappa} V_{\nu}^*.\]
 The following von Neumann type inequality is now an immediate consequence of Theorem \ref{tPoisson}.
Note that it implies in particular that the definition of $p (\Vcont)$ does not depend on the representation chosen in
\eqref{poly}.

\begin{cor} \label{vNi}
Let $\La$ be finitely aligned, $\Vcont$ be a $\La$-contraction satisfying the Popescu condition and let
$p$ be a polynomial in noncommuting variables indexed by $\La \times \La$. Then the following von Neumann
type inequality holds: \[ \|p(\Vcont)\| \leq \|p({\mathcal L})\|.\]
\end{cor}

When $\La$ is finite with no sources and $\Vcont$ is a $\La$-isometry, the Poisson transform associated to $\Vcont$
may be regarded as a map on $\Cun$:

\begin{tw} \label{cPoisson}
Let $\La$ be finite and without sources and let $\Vcont$ be a $\La$-isometry. Then there exists a unique unital completely positive map $T_{\Vcont}:\Cun
\to B(\Hil)$ satisfying
\[T_{\Vcont}(s_{\la} s_{\mu}^*) = V_{\la} V_{\mu}^*, \;\;\;\; \la, \mu \in \La.\]
\end{tw}
\begin{proof}
It is enough to check that the map $\Pois$ constructed in the previous theorem vanishes on the ideal $\mathcal{J}$. To this end
consider an operator of the form $xP_j y$, where $x,y \in \Toep$, $j\in \oner$. We can assume that $x=L_{\lambda} L_{\nu}^*$,
$y=L_{\alpha} L_{\beta}^*$ ($\lambda, \nu, \alpha, \beta \in \La$). Further we can also assume that in fact $xP_j y= L_{\mu}P_j
L_{\nu}^*$: an operator $P_j L_{\alpha}$ is non-zero only if $\sigma(\alpha)_j=0$, and then $P_j L_{\alpha} = L_{\alpha} P_j$.
But then, as $P_j = I - \sum_{\gamma \in \La^{e_j}} L_{\gamma}L_{\gamma}^*$, there is
\begin{align*} \Pois (xP_j y) &= \Pois ( L_{\mu} L_{\nu}^* - \sum_{\gamma
\in \La^{e_j}} L_{\mu\gamma}L_{\nu\gamma}^* ) = V_{\mu} V_{\nu}^* - \sum_{\gamma \in \La^{e_j}}
V_{\mu\gamma}V_{\nu\gamma}^* \\&= V_{\mu} V_{\nu}^* - V_{\mu} \left(\sum_{\gamma \in \La^{e_j}}
V_{\gamma}V_{\gamma}^* \right) V_{\nu}^* = 0.\end{align*}
\end{proof}

The following corollary will be of use for the analysis of states on $\Cun$ in Section 3.

\begin{cor}   \label{corPois}
Let $\La$ be finite without sources and let $\Vcont=\{V_{\la}: \la \in \La\}$ be a family of operators satisfying conditions (i)-(ii) of Definition \ref{Lcont}. Let $D\in B(\Hil)$ be a positive operator such that
\begin{equation} \forall_{n\in \Nr} \;\; \sum_{\la \in \La^n} V_{\la} D V_{\la}^* = D.
\label{Dinv}\end{equation}
Then there exists a unique completely positive map $T_{\Vcont,D}:\Cun \to B(\Hil)$ satisfying
\[T_{\Vcont,D}(s_{\la} s_{\mu}^*) = V_{\la} D V_{\mu}^*, \;\;\;\; \la, \mu \in \La.\]
\end{cor}

\begin{proof}
Following \cite{pureod} in the rank 1 case assume for the moment that $D$ has a
bounded inverse and consider the $\La$-isometry given by the family
$\{D^{-\frac{1}{2}} V_{\la} D^{\frac{1}{2}} : \la \in \La\}$. The construction from the
previous theorem yields the unital completely positive map $T: \Cun \to B(\Hil)$ such that
\[T(s_{\la} s_{\mu}^*) = D^{-\frac{1}{2}}  V_{\la} D V_{\mu}^* D^{-\frac{1}{2}} ,
\;\;\;\; \la, \mu \in \La.\] It is easy to see that the map defined by
\[ T_{\Vcont,D} (x) = D^{\frac{1}{2}} T(x)  D^{\frac{1}{2}}, \;\;\;
x \in \Cun,\] satisfies the condition in the corollary.

In the general case we may modify (following \cite{ergod} in the rank 1 case) the defect
operator and Poisson transform as follows. Define
\begin{align*}
\Delta_{s,D}(\Vcont) &= \sum_{\mu \in \La, \, \sigma(\mu) \leq e} (-s^2)^{|\mu|} V_{\mu} D V_{\mu}^*\\
&= \sum_{n \leq e} (-s^2)^{|n|} \sum_{\mu \in \La^n} V_{\mu} D V_{\mu}^*\\
&= \sum_{n \leq e} (-s^2)^{|n|} D = (1-s^2) D
\end{align*}
which is positive for all $s \in (0,1)$. The same arguments as before show that
\[
\sum_{\la \in \La} s^{2|\la|} V_{\la} \Delta_{s,D}(\Vcont)  V_{\la}^*  = D
\]
so that $W_{s,D} : \Hil \to \ell^2(\La) \otimes \Hil$ given by $W_{s,D}(\xi)= \sum_{\la \in \La} \delta_{\la} \ot s^{|\la|} \Delta_{s,D}^{\frac{1}{2}} V_{\la}^* \xi$ verifies $W_{s,D}^* W_{s,D}=D$ and
$T_{\Vcont,D}(x)= \lim_{s \to 1}W_{s,D}^* (x \otimes I)W_{s,D}$ is the required completely positive map.
%it is enough to take $\epsilon>0$ and replace (for the moment) the operator $D$ with
%(boundedly invertible) $D_{\epsilon} =D + \epsilon I_K$. The maps $T_{\Vcont,D_{\epsilon}}$ may be shown to
%converge pointwise as $\epsilon$ tends to $0^+$; the limit yields the required map.
\end{proof}

Note that \eqref{Dinv} may be interpreted as the statement that the operator $D$ is left invariant
by certain natural completely positive maps acting on $B(\Hil)$. This will be exploited later.

\section{Dilation and commutant lifting theorems for $\La$-contractions}

In this section we construct (minimal) dilations of $\La$-contractions satisfying the Popescu condition
and of $\La$-isometries. Related commutant lifting theorems are also discussed and canonical unital completely
positive maps associated to a given $\La$-isometry introduced.

\subsection*{Dilations}

The Poisson transform in conjunction with the Stinespring Theorem provides a dilation of $\La$-contractions
(or $\La$-isometries) to Toeplitz-Cuntz-Krieger (or Cuntz-Pimsner) families. This can be thought of as a dilation
of a family of contractions satisfying certain commutation relations to a family of partial isometries which will
not only satisfy the initial commutation relations but also extra conditions involving their adjoints. As is well
known, although the dilation of a triple of commuting contractions to commuting isometries may not be possible, it may
always be constructed as long as the contractions in question are doubly commuting. This is reflected in our context
by the requirement that the $\La$-contractions we are dilating are supposed to satisfy the Popescu condition.

\begin{tw} \label{Tdilation}
Let $\La$ be finitely aligned and let $\Vcont$ be a $\La$-contraction on a Hilbert space $\Hil$
satisfying the Popescu condition. There exists a
Hilbert space $\Kil \supset \Hil$ and a $\La$-contraction $\Wcont$ on $\Kil$
consisting of partial isometries forming a Toeplitz-Cuntz-Krieger family such that for each $\la \in \La$
\[ W_{\la}^*|_{\Hil} = V_{\la}^*.\]
One may assume that $\Kil = \overline{\textup{Lin}} \{W_{\la} \Hil: \la \in \La\}$; under this assumption
the family $\Wcont$ is unique up to unitary equivalence and is called the minimal dilation of
$\Vcont$ (to a Toeplitz-Cuntz-Krieger family).
\end{tw}
\begin{proof}
Consider the minimal Stinespring dilation of the Poisson transform $\Pois$ constructed in Theorem \ref{tPoisson}.
This provides us with a Hilbert space $\Kil$, a representation $\pi: \Toep \to B(\Kil)$ and an operator
$V\in B(\Hil;\Kil)$ such that for all $x \in \Toep$
\[ \Pois(x) = V^* \pi(x) V\]
and $\Kil = \overline{\Lin} \{\pi(x)V\xi: x \in \Toep, \xi \in \Hil\}$. We may assume that $V$ is an isometry,
as even when $\La^0$ is infinite we can `make' $\Pois$ and $\pi$ unital by passing to $\Toep^1$. This allows us to view
$\Hil$ as a subspace of $\Kil$.
Define for each $\la \in \La$
\[ W_{\la} = \pi (L_{\la}).\]
It is clear that the family $\{W_{\la}: \la \in \La\}$ is a Toeplitz-Cuntz-Krieger family, so in particular a $\La$-contraction.
Let $\la, \mu, \nu \in \La$ and $\xi, \eta \in \Hil$. Then
\begin{align*} \langle W_{\la}^* \xi, \pi(L_{\mu} L_{\nu}^*) \eta \ra &=
 \langle \xi, W_{\la} W_{\mu} W_{\nu}^* \eta \ra =
\langle \xi, W_{\la\mu} W_{\nu}^* \eta \ra \\&= \langle \xi, V_{\la\mu} V_{\nu}^* \eta \ra =
\langle V_{\la}^* \xi, V_{\mu} V_{\nu}^* \eta \ra =
  \langle V_{\la}^* \xi, \pi(L_{\mu} L_{\nu}^*) \eta \ra \end{align*}
and by minimality this implies that each $W_{\la}^*$ leaves $\Hil$ invariant and $W_{\la}^*|_{\Hil} = V_{\la}^*$.
The minimality condition may be therefore written as
$\Kil = \overline{\Lin} \{W_{\la}\xi: \la \in \La, \xi \in \Hil\}$.
The uniqueness claim follows since one can easily check that if $\Wcont'$ is a $\La$-contraction on a Hilbert space $\Kil'$
satisfying the requirements in the theorem, then the continuous linear extension of the map $W_{\la} \xi \mapsto
W'_{\la} \xi$ ($\la \in \La, \xi \in \Hil$) yields the intertwining unitary from $\Kil$ to $\Kil'$.
\end{proof}

Note that this dilation result in particular subsumes the dilation theorem for contractive $A$-relation tuples
obtained in \cite{RajJSan}.
%When the graph $\La$ has only one vertex, it reduces to a special case of Theorem %4.7 of\cite{Solel}.
It is also closely connected to general dilation theory for completely contractive representations of product systems of $C*$-correspondences
(\cite{Solel}). For the description of these connections we refer to the
forthcoming note \cite{Adam}, where it is in particular established that
every rank-2 contraction can be dilated to a Toeplitz type family.

If $\Vcont$ is a $\La$-isometry it is natural to expect that it can be dilated to a Cuntz-Pimsner $\La$-family.
The theorem below shows that this is automatically verified for the minimal dilation to
a Toeplitz-Cuntz-Krieger $\La$-family.

\begin{tw} \label{mintoepcun}
Assume that $\La$ is row-finite. Then the minimal dilation of a $\La$-isometry is a Cuntz-Pimsner $\La$-family.
%on a Hilbert space $\Kil$.
%satisfying the Popescu condition.
%Suppose that $\Hil$ is a Hilbert space and $\Wcont$ is a $\La$-contraction  on $\Hil$
%being the minimal dilation of $\Vcont$ to a Toeplitz-Cuntz-Krieger family.
%Then $\Wcont$ is a $\La$-isometry and a Cuntz-Pimsner $\La$-family.
\end{tw}

\begin{proof}
Let $\Vcont$ be a $\La$-isometry on $\Hil$.
As row-finite higher rank graphs are finitely aligned and $\La$-isometries
satisfy the Popescu condition, we can apply Theorem \ref{Tdilation} to $\Vcont$ to obtain its minimal dilation
$\Wcont$ on a Hilbert space $\Kil \supset \Hil$. By minimality it suffices to show the following:
\begin{equation} \langle W_{\gamma} \xi, W_a W_{\nu} \eta \ra = \sum_{\la \in \La_a^n}
\langle W_{\la} W_{\la}^* W_{\gamma} \xi, W_a W_{\nu} \eta \ra \label{long}\end{equation}
for all $\gamma, \nu \in \La$, $a \in \La^0, n \in \Nr$ and $\xi, \eta \in \Hil$. We can assume that
$r(\nu)=a$. Moreover, using the fact that $\Wcont$ is a Toeplitz-Cuntz-Krieger family,
we see that $W_a W_{\la} W_{\la}^* = 0 $, whenever $\la \in \La$ and $r(\la) \neq a$ (using \ref{TCK} (i) and (iii)).
This, together with \ref{TCK} (v) implies that the right hand side of \eqref{long} is  equal to
\begin{align*} \sum_{\la \in \La^n}
\langle W_{\la}^* W_{\gamma} \xi  ,  W_{\la}^*  W_{\nu} & \eta \ra
= \sum_{\la \in \La^n} \;\;\sum_{\la \alpha = \gamma \beta \in MCE(\la, \gamma)}
\langle W_{\alpha} W_{\beta}^*\xi, W_{\la}^*  W_{\nu} \eta \ra \\ =&
\sum_{\la \in \La^n} \;\; \sum_{\la \alpha = \gamma \beta \in MCE(\la, \gamma)}
\langle  W_{\beta}^*\xi, W_{\la \alpha}^*  W_{\nu} \eta \ra \\ =&
\sum_{\la \in \La^n} \;\; \sum_{\la \alpha = \gamma \beta \in MCE(\la, \gamma)}
\;\; \sum_{\la \alpha \mu= \nu \kappa \in MCE(\la\alpha, \nu)}
\langle  W_{\beta}^*\xi, W_{\mu}  W_{\kappa}^* \eta \ra \\=&
\sum_{\la \in \La^n} \;\; \sum_{\la \alpha = \gamma \beta \in MCE(\la, \gamma)}
\;\;\sum_{\la \alpha \mu= \nu \kappa \in MCE(\la\alpha, \nu)}
\langle  W_{\beta \mu}^*\xi,   W_{\kappa}^* \eta \ra \\ =&
\sum_{\la \in \La^n} \;\; \sum_{\la \alpha = \gamma \beta \in MCE(\la, \gamma)}
\;\;\sum_{\la \alpha \mu= \nu \kappa \in MCE(\la\alpha, \nu)}
\langle  V_{\beta \mu}^*\xi,  V_{\kappa}^* \eta \ra.
\end{align*}
Suppose that $\sigma(\gamma) = m$, $\sigma(\nu)=p$. Then the factorisation property implies that
the sum above can be rewritten as follows:
\begin{align*} \sum_{\la \in \La^n}
\langle W_{\la}^* W_{\gamma} \xi  ,  W_{\la}^*  W_{\nu} & \eta \ra
= \sum_{\beta \in \La^{n \vee m - m }, r(\beta) =s(\gamma)}
\;\;\sum_{\gamma \beta \mu= \nu \kappa \in MCE(\gamma \beta, \nu)}
\langle  V_{\beta \mu}^*\xi,  V_{\kappa}^* \eta \ra
\\ =& \sum_{\beta \in \La^{n \vee m - m }, r(\beta) =s(\gamma)}
\;\;\sum_{\gamma \delta = \nu \kappa \in \La^{n\vee m \vee p}}
\langle  V_{\delta}^*\xi,  V_{\kappa}^* \eta \ra
%\end{align*}
\\=& \sum_{\gamma \alpha' = \nu \beta' \in MCE(\gamma, \nu)}
\;\; \sum_{\kappa'\in \La^{n\vee m \vee p - n\vee m}}\langle  V_{\kappa'}^*V_{\alpha'}^*\xi,  V_{\kappa'}^* V_{\beta'}^* \eta \ra
\\=& \sum_{\gamma \alpha' = \nu \beta' \in MCE(\gamma, \nu)}
\;\;  \left \langle\sum_{\kappa'\in \La^{n\vee m \vee p - n\vee m}}  V_{\kappa'} V_{\kappa'}^*V_{\alpha'}^*\xi \; , \;  V_{\beta'}^*\eta \right\ra.
\end{align*}
As $\Vcont$ is a $\La$-isometry we finally obtain
\[ \sum_{\la \in \La^n}
\langle W_{\la}^* W_{\gamma} \xi  ,  W_{\la}^*  W_{\nu} \eta \ra
= \sum_{\gamma \alpha' = \nu \beta' \in MCE(\gamma, \nu)}
\langle  V_{\alpha'}^*\xi, V_{\beta'}^*\eta \ra. \]
Using again condition (v) in the definition of a Toeplitz-Cuntz-Krieger family it is easy to check that the sum above is equal to the expression on the right hand side of \eqref{long}.

As $\Wcont$ is a Cuntz-Pimsner $\La$-family, to show that it is a $\La$-isometry it is enough to prove that
$\sum_{a \in \La^0} W_a = I_{\Hil}$. This follows immediately from the minimality condition, as
for all $\gamma\in \La$, and $\xi\in \Hil$
\[ \sum_{a \in \La^0} W_a W_{\gamma} \xi = W_{r(\gamma)} W_{\gamma} \xi = W_{\gamma} \xi.\]
\end{proof}

A similar result for rank-2 graphs with one vertex has been established in Lemma 5.2 of \cite{2dilat}. Theorem \ref{mintoepcun}
 implies the following corollary.

\begin{cor} \label{Cdilation}
Let $\La$ be finite and let $\Vcont$ be a $\La$-isometry on a Hilbert space $\Hil$.
There exists a
Hilbert space $\Kil \supset \Hil$ and a $\La$-isometry $\Wcont$ on $\Kil$
consisting of partial isometries forming a Cuntz-Pimsner family such that for each $\la \in \La$
\[ W_{\la}^*|_{\Hil} = V_{\la}^*.\]
One may assume that $\Kil = \overline{\textup{Lin}} \{W_{\la} \Hil: \la \in \La\}$; under this assumption
the family $\Wcont$ is unique up to unitary equivalence.
\end{cor}

Note that if $\La$ is additionally assumed to have no sources the corollary can be proved directly along identical lines as
Theorem \ref{Tdilation}, this time exploiting the version of the Poisson transform obtained in Theorem \ref{cPoisson}.

%It actually turns out that Theorem \ref{Cdilation} may be viewed as a special case of
%Theorem \ref{Tdilation}. This is formalised in the next proposition.

\subsection*{Commutant lifting theorems}

Our first commutant lifting theorem concerns dilations of $\La$-contractions. It is an immediate  consequence
of Arveson's  commutant lifting result (Theorem 1.3.1 of \cite{subalg})
exploiting the way in which the dilations were constructed (see also \cite{PPois}).
Whenever $\Vcont$ is a $\La$-contraction on a Hilbert space $\Hil$ we will write
$\Vcont' = \{ T \in B(\Hil): \forall_{\la \in \La} T V_{\la} = V_{\la}T \}$. $\Vcont'$ is generally a non-selfadjoint operator algebra,
whereas $(\Vcont \cup \Vcont^*)'=\Vcont' \cap (\Vcont^*)'$ is a von Neumann algebra, further by a slight abuse of language called the \emph{commutant} of $\Vcont$ .

\begin{tw} \label{cl1}
Let $\La$ be finitely aligned and let  $\Vcont$ be a $\La$-contraction on $\Hil$
satisfying the Popescu condition. Let $\Wcont$ be its minimal dilation to a
Toeplitz-Cuntz-Krieger family acting on a Hilbert space $\Kil \supset \Hil$ and let $P\in B(\Kil)$ denote the orthogonal
projection onto $\Hil$. For any $X \in \Vcont'\cap (\Vcont^*)' \subset B(\Hil)$ there exists a unique
$\wt{X} \in (\Wcont \cup \Wcont^*)' \cap \{P\}'$ such that $ \wt{X}|_{\Hil}= X$. The correspondence $X \mapsto \wt{X}$ is
a normal $*$-isomorphism of the von Neumann algebras $(\Vcont \cup \Vcont^*)'$ and $(\Wcont \cup \Wcont^* \cup \{P\})'$.
\end{tw}

Before we formulate the second of our commutant lifting theorems we need to discuss
natural families of unital completely positive maps arising from $\La$-isometries,
which generalise familiar endomorphisms of $B(\Hil)$ associated with representations
of Cuntz algebras on a Hilbert space $\Hil$.

  \begin{deft} \label{defcp}
Suppose that $\La$ is cofinal  and  $\Vcont$ is a $\La$-contraction on a Hilbert space $\Hil$.
Define for each $n \in \Nr$ and $ X \in B(\Hil)$
\[ \sigma_{\Vcont} (n) (X) = \sum_{\la \in \La^n} V_{\la} X V_{\la}^*\]
(cofinality assures that each $\La^n$ is non-empty).
Each $\sigma_{\Vcont}(n)$ is a completely positive contraction on $B(\Hil)$ and moreover for all
$n, m \in \Nr$
\[ \sigma_{\Vcont} (n+m) = \sigma_{\Vcont} (n) \circ \sigma_{\Vcont} (m).\]
The resulting action of $\Nr$ on $B(\Hil)$ will be denoted by $\sigma_{\Vcont}$.
It is unital if and only if $\Vcont$ is a $\La$-isometry.
We also introduce a selfadjoint subspace of operators left invariant by the action:
\[ \Fix \sigma_{\Vcont}= \{X \in B(\Hil):
 \forall_{n \in \Nr}\;\;  \sigma_{\Vcont} (n) (X) = X\}.\]
%\[ \Fix \sigma_{\Vcont}= \{X \in B(\Hil):
% \forall_{m \in \Nr}\;\;  \sigma_{\Vcont} (m) (X) = X\}.\]
\end{deft}

%To obtain interesting objects we will have to assume generally that $\La^n \neq \emptyset$ for at least some
%or better all non-zero $n$. To require $\La$ to be cofinal seems a reasonable assumption. In that case it is clear that
%$\sigma_{\Vcont}(n)$ is non-zero for all $n$.

Note that $\Fix \sigma_{\Vcont} \subset \bigcap_{j=1}^r \Fix \sigma_{\Vcont}(e_j)$ and an operator
in $\bigcap_{j=1}^r \Fix \sigma_{\Vcont}(e_j)$ belongs to $\Fix \sigma_{\Vcont}$ if and only if it is
`diagonal' with respect to the decomposition $\Hil = \bigoplus_{a \in \La^0} V_a \Hil$. In fact
the inclusion above is an equality except in degenerate cases. As it is often implicitly used in what follows we formulate it as a lemma.

%In general these operator systems do not have to be equal but this can only happen for
%$\La$-contractions.
%is not sufficiently connected: if a
%rank-one graph $\La$  has at least 2 vertices and no connecting edges,
%and $\Vcont$ is a nontrivial $\La$- on a Hilbert space $\Hil$, then
%$\Fix \sigma_{\Vcont}=B(\Hil) \neq \Fizx \sigma_{\Vcont}$.
%The following lemma shows that the only obstructions for $\Fix \sigma_{\Vcont}= \Fizx \sigma_{\Vcont}$
%may essentially arise from the situation as in the example above
%(i.e.\.insufficient connectivity between vertices).

\begin{lem} \label{same}
Suppose that $\La$ is cofinal and $\Vcont$ is a $\La$-contraction on $\Hil$. Then
$\textup{Fix} \: \sigma_{\Vcont} = \bigcap_{j=1}^r \textup{Fix} \: \sigma_{\Vcont}(e_j)$.
%\textup{Fix} \sigma_{\Vcont}$.
\end{lem}

\begin{proof}
Suppose  that  $X \in \Fix \sigma_{\Vcont}(e_j)$ for all $ j \in \oner$ so that
$X = \sum_{\la \in \La^{e_j}} V_{\la} X V_{\la}^*.$
Thus we have
\[
\sum_{a \in \La^0} V_a XV_a = \sum_{\la \in \La^{e_j}, \: a \in \La^0} V_a V_{\la} X (V_a V_{\la})^* =
\sum_{\la \in \La^{e_j}} V_{\la} X V_{\la}^*=X,
\]
so that $X= \sum_{a \in \La^0} V_a XV_a = \sigma_{\Vcont}(0) (X)$ and in particular $X$ is diagonal with respect to the decomposition
$\Hil = \bigoplus_{a \in \La^0} V_a \Hil$. Moreover, $\sigma_{\Vcont}(n) (X)= \sigma_{\Vcont}(n_1 e_1) \circ  \ldots \circ \sigma_{\Vcont}(n_r e_r)(X)= X$ for all $n \in \Nr$.
\end{proof}

%and multiplying the equality above by $V_{\nu}$ (where $\nu \in \La^{e_j}$) on the right we obtain
%\begin{equation} \label{com1} X V_{\nu} = \sum_{\la \in \La^{e_j}} V_{\la} X V_{\la}^* V_{\nu}.\end{equation}
%Let now $a, b \in \La^0$, $a \neq b$. As $\La$ has no sources, there is $\nu \in \La^{e_j}$ such that
%$a = r(\nu)$. For such $\nu$
%\[  V_b X V_a V_{\nu} =  V_b X V_{\nu} = V_b \sum_{\la \in \La^{e_j}} V_{\la} X V_{\la}^* V_{\nu}  =
%  \sum_{\la \in \La^{e_j}_b} V_b V_{\la} X V_{\la}^* V_b V_{r(\nu)}V_{\nu}0 .\]
%If now $\nu\in \La^{e_j}$ and $r(\nu) \neq a $ then we also have $V_b X V_a V_{\nu} = 0$,
%and this implies that
%\[ V_b X V_a = \sum_{\nu \in \La^{e_j}} V_b X V_a V_{\nu} V_{\nu}^* = 0,\]
%so that $X \in \Fizx \sigma_{\Vcont}$.

The next lemma generalises a well known result on connections between the space of fixed points of an endomorphism of $B(\Hil)$
and the commutant of the corresponding representation of a Cuntz algebra ${\mathcal O}_n$.

\begin{lem}       \label{fix1}
Suppose that $\La$ is  finite  and has no sources. Let $(\pi, \Kil)$ be a representation of $\Cun$ and
let $\Vcont$ be the $\La$-isometry given by $V_{\la} = \pi (s_{\la})$.
Then $\Fix \, \sigma_{\Vcont}  = \pi (\Cun)'.$
\end{lem}

\begin{proof} %By \ref{same} we know that $\Fix \sigma_{\Vcont} = \Fizx \sigma_{\Vcont}$.
Let $X \in \Fix \sigma_{\Vcont}$ and $\nu \in \La$. Then $XV_{\nu} = \sum_{\la \in \La^{\sigma(\nu)}} V_{\la} X V_{\la}^* V_{\nu} = V_{\nu} X$
%Let $j \in \oner$ and  $\nu \in \La^{e_j}$. Arguing as in the last
%lemma we deduce that formula \eqref{com1} holds and $X$ is `diagonal'. Combining these facts together we obtain
%\[ X V_{\nu} = V_{\nu} X V_{s(\nu)} = V_{\nu} V_{s(\nu)}  X V_{s(\nu)} =
%V_{\nu} \sum_{a\in \La^0} V_a X V_a = V_{\nu}X\]
%for all $\nu \in \La^{e_j}$. As $j \in \oner$ was arbitrary, we obtain $XV_{\la} = V_{\la}  X$
%for all $\la \in \La$.
Moreover as $\Fix \sigma_{\Vcont}$ is selfadjoint, we also get
$XV_{\nu}^* = V_{\nu}^*  X$, hence $X \in \pi(\Cun)'$.

The inclusion $\pi (\Cun)' \subset \Fix \sigma_{\Vcont}$ is obvious.
\end{proof}

%Note that the proofs of both lemmas above could be conducted for each `direction' separately. It is essentially enough
%to establish the result for rank-one graphs.

The following lemma is an extension of the well known result concerning states and GNS representations
(see for example Proposition 3.10 in \cite{Takes}).
It follows from results in \cite{subalg} Section 1.4. We include a short proof for the readers convenience.

\begin{lem}   \label{funlem}
Let $\alg$ be a unital $C^*$-algebra, $\Hil$ be a Hilbert space and suppose we are given two completely positive maps
$\Psi, \Phi: \alg \to B(\Hil)$. Assume that $\Phi$ is unital, $\Phi-\Psi$ is completely positive and
let $\Phi=V^*\pi(\cdot) V$ be the minimal Stinespring decomposition of $\Phi$ (so that $(\pi, \Kil)$ is a representation
of $\alg$, $V: \Hil \to \Kil$ is an isometry and $\Kil = \ol{\textup{Lin}}\{\pi(a) V \xi: a \in \alg, \xi \in \Hil\}$).
Then there exists a unique operator $X \in \pi(\alg)'\subset B(\Kil)$ such that
\begin{equation} \Psi(a) = V^* (X \pi(a)) V, \;\;\; \text{for all } a \in \alg.
\label{funnylem}\end{equation}
Moreover $0\leq X \leq I$.
\end{lem}

\begin{proof}
Write $\pi(\alg) \Hil$ for $\Lin\{\pi(a) V \xi: a \in \alg, \xi \in \Hil\}$.
Consider a quadratic form on $ \pi(\alg) \Hil$ given by
\[ B\left(\sum_{i=1}^k \pi(a_i) \xi_i \right) = \sum_{i,j=1}^k \langle \xi_i , \Psi(a_i^* a_j) \xi_j \rangle,\]
$k \in \bn$, $a_1, \ldots, a_k \in \alg$, $\xi_1, \ldots \xi_k \in \Hil$.
Note that as $\Psi$ is completely positive, $B$ is positive-definite. Moreover
as $\Phi-\Psi$ is completely positive
\[ \left\|\sum_{i=1}^k \pi(a_i) \xi_i \right\|^2 -   B\left(\sum_{i=1}^k \pi(a_i) \xi_i \right) =
 \sum_{i,j=1}^k \langle \xi_i , (\Phi-\Psi)(a_i^* a_j) \xi_j \rangle \geq 0,\]
so that  the sesquilinear form $B': \pi(\alg) \Hil \times\pi(\alg) \Hil$
given by
\[ B'\left( \sum_{i=1}^k \pi(a_i) \xi_i, \sum_{j=1}^k \pi(b_j) \eta_j \right) =
\sum_{i,j=1}^k \langle \xi_i , \Psi(a_i^* b_j) \eta_j \rangle,\]
($k \in \bn$, $a_1, \ldots, a_k,b_1, \ldots, b_k \in \alg$,
$\xi_1, \ldots \xi_k, \eta_1, \ldots \eta_k\in \Hil$) represents a bounded operator
on $\Kil$. More precisely, there exists a unique $X \in B(\Kil)$ such that
for all $\zeta_1, \zeta_2 \in \pi(\alg) \Hil$
\[\langle \zeta_1, X \zeta_2 \ra = B'(\zeta_1, \zeta_2).\]
It is now elementary to check that $X$ is a positive contraction in $\pi(\alg)'$
satisfying \eqref{funnylem}.
\end{proof}

Note that it is not enough to assume in the lemma above that $\Phi-\Psi$ is positive; it is possible to find examples of
completely positive maps $\Phi$, $\Psi$ such that $\Phi-\Psi$ is positive but not completely positive. In such case it is clearly
impossible to have a representation as in \eqref{funnylem} for a positive $X\in \pi(\alg)'$.

We are now ready to state and prove
the second of the commutant lifting theorems in this section. It concerns
dilations of $\La$-isometries (as opposed to $\La$-contractions). When $\Vcont$ is a $\La$-isometry
on $\Hil$, Theorem \ref{cl1} implies that the commutant of $\Vcont$ is isomorphic to the intersection of the commutant
of the representation of $\Cun$ induced by the minimal dilation of $\Vcont$ with the algebra of operators
diagonal with respect to the decomposition of the dilation space: $\Kil=\Hil \oplus (\Kil \ominus \Hil)$.
The theorem below shows that we can
actually find an alternative representation of  the selfadjoint part of the whole commutant of the afore-mentioned representation of $\Cun$. The way to do it is suggested by the Lemma \ref{fix1}, which shows that there is a close connection between the commutants we are interested in and fixed points of the relevant completely positive maps.
As the fixed point subspaces of completely positive maps do not have to be closed under multiplication, it is natural that here the identification obtained may be valid only in the category of operator systems (i.e.\ a $^*$-homomorphism is replaced by an isometric order isomorphism). For more discussion on this topic
we refer the reader to the paper \cite{pureod}.

\begin{tw}
Suppose that $\La$ is finite without sources. Let $\Vcont$ be a $\La$-isometry on a Hilbert space $\Hil$ and let $\Wcont$ be its
minimal dilation to a Cuntz-Pimsner family on $\Kil$. Let $P:\Kil \to \Hil$ denote the projection onto $\Hil$ and let $(\pi,
\Kil)$ be the representation of $\Cun$ determined by $\Wcont$. Then the map $X \mapsto PXP$ yields a  complete order isomorphism
between the  commutant $\pi(\Cun)'$ and the  operator system $\Fix \sigma_{\Vcont}$.
%(complete order isomorphism here is understood as a fact that both the map
%and its inverse are completely positive).
\end{tw}

\begin{proof}
Suppose that $X \in \pi(\Cun)'$. To show that $PXP \in \Fix \sigma_{\Vcont}$ notice that by Corollary \ref{Cdilation} we have
$PW_{\la} = PV_{\la}P= V_{\la}P$ for all $\la \in \La$ hence
\[
\sigma_{\Vcont}(n)(PXP) = \sum_{\la \in \La^n} V_{\la} PXP V_{\la}^* = \sum_{\la \in \La^n} PW_{\la} X W_{\la}^*P =
\sum_{\la \in \La^n} PW_{\la} W_{\la}^*XP=PXP
\]
for $n \in \Nr$.

% \Fizx \sigma_{\Vcont}$
%it suffices to show that for arbitrary  $\nu, \mu \in \La$, $j \in \oner$ and
%$\xi, \eta \in \Hil$,
%\begin{equation}
%\langle W_{\nu}^* \xi, \sigma_{\Vcont}(e_j) (PXP) W_{\mu}^* \eta \ra =
% \langle W_{\nu}^* \xi,  X W_{\mu}^* \eta \ra.
%\label{PXP} \end{equation}
%As we have
%\begin{align*}
%\langle W_{\nu}^* \xi, & \sigma_{\Vcont}(e_j) ( PXP) W_{\mu}^* \eta \ra =
%\langle W_{\nu}^* \xi, \sum_{\la \in \La^{e_j}}  W_{\la} PXP W_{\la}^* W_{\mu}^* \eta \ra
%\\ =&  \sum_{\la \in \La^{e_j}} \langle W_{\nu \la}^* \xi,   PXP  W_{\mu\la}^* \eta \ra=
% \sum_{\la \in \La^{e_j}} \langle W_{\nu \la}^* \xi,   X  W_{\mu\la}^* \eta \ra
%\\ =&  \sum_{\la \in \La^{e_j}} \langle W_{\nu}^* \xi, W_{\la} W_{\la}^*  X  W_{\mu}^* \eta \ra
%= \langle W_{\nu}^* \xi, X  W_{\mu}^* \eta \ra,
%\end{align*}
%so that \eqref{PXP} is proved.

%Obviously if $X$ is selfadjoint (or positive), so is $PXP$.

Conversely suppose that $D \in \Fix \sigma_{\Vcont}$, and assume first $0 \leq D \leq I_{\Hil}$.
Corollary \ref{corPois}
implies the existence of completely positive maps $
T_{\Vcont,I-D},T_{\Vcont,D}: \Cun \to B(\Hil)$
such that
\[T_{\Vcont,D}(s_{\la} s_{\mu}^*) = V_{\la} D V_{\mu}^*, \;\; T_{\Vcont,I-D}(s_{\la} s_{\mu}^*)
= V_{\la} (I-D) V_{\mu}^* \;\; \la, \mu \in \La. \]
It is easy to see that $T_{\Vcont,I-D} = T_{\Vcont} - T_{\Vcont,D}$.
Recall that the minimal dilation of $\Vcont$ to a $\La$-isometry was achieved via the
Stinespring dilation
of the map $T_{\Vcont}$. By Lemma \ref{funlem}  there
exists a unique operator $X \in \pi(\Cun)'$ such that  $0 \leq X \leq I_{\Kil}$   and
\[ \langle \xi, T_{\Vcont,D}(S) \eta \ra = \langle \xi , X \pi(S) \eta \rangle, \;\;\;
S \in \Cun, \xi, \eta \in \Hil.\]
Comparing the formulas above yields the equality
\[ \langle \xi, V_{\la} D V_{\mu}^* \eta \ra =  \langle \xi, W_{\la} X W_{\mu}^* \eta \ra,\]
for arbitrary $\la, \mu \in \La$ and $\xi, \eta \in \Hil$, so that $D=PXP$.

The fact that the correspondence $X \mapsto PXP$ preserves order and norm when
restricted to respective selfadjoint parts may be now established exactly as in
Proposition 4.1 of \cite{pureod}.

A unital isometric order isomorphism $\phi_h$ between selfadjoint parts of two operator systems $Y_1$ and $Y_2$ has a unique extension to a complex linear map $\phi:Y_1 \to Y_2$. It is easy to check that $\phi$ is a (Banach space) isomorphism; moreover, as it is positive and unital, it has to be contractive. The same argument applied to the inverse of $\phi$ shows that actually $\phi$ has to be isometric.

It remains to show that all the above properties of the map $P \mapsto PXP$ remain valid when we pass to its matrix liftings mapping $M_n(\pi(\Cun)') \to M_n(\Fix \sigma_{\Vcont})$ ($n \in \bn$).

To this end fix $n \in \bn$ and consider the $\La$-contraction $\Vcont^{(n)}$ on $\Hil^{\oplus n}$ defined by $V_{\la}^{(n)}=
V_{\la} \oplus \cdots \oplus V_{\la}$, ($\la \in \La$). The minimal dilation of $\Vcont^{(n)}$ to a Cuntz-Pimsner family can be
identified with the $\La$-contraction $\Wcont^{(n)}$. Let $P^{(n)}: \Kil^{\oplus n} \to \Hil^{\oplus n}$ be the relevant
orthogonal projection and let $\pi^{(n)}$ denote the representation of $\Cun$ on $\Kil^{\oplus n}$ associated with
$\Wcont^{(n)}$. From the first part of the proof it follows that the map $X^{(n)} \mapsto P^{(n)} X^{(n)} P^{(n)}$ is an
isometric unital order isomorphism of $(\pi^{(n)}(\Cun))'$ and $\Fix \sigma_{\Vcont^{(n)}}$. It remains to note that actually
$(\pi^{(n)}(\Cun))' = (I_n \otimes \pi(\Cun))' = M_n \otimes \pi(\Cun)'$, $\Fix \sigma_{\Vcont^{(n)}} = M_n (\Fix
\sigma_{\Vcont})$ and the map $X^{(n)} \mapsto P^{(n)} X^{(n)} P^{(n)}$ is equal to the matrix lifting of $P \mapsto PXP$. This
ends the proof.
\end{proof}

The theorem above is a generalisation of Theorem 5.1 in \cite{pureod}. Using the methods identical
to the ones of that paper we can establish the following corollary.

\begin{cor}
Let $\La$ be finite without sources. Suppose that  $\Vcont$ and $\wt{\Vcont}$ are $\La$-isometries on a Hilbert space $\Hil$ and
let $\Wcont$, $\wt{\Wcont}$ be their respective minimal dilations to Cuntz-Pimsner families (acting respectively on Hilbert
spaces $\Kil$ and $\wt{\Kil}$). There is a completely isometric correspondence between the set of intertwiners between $\Wcont$
and $\wt{\Wcont}$ (i.e.\ operators $U \in B(\Kil; \wt{\Kil})$ such that for each $\la \in \La$ there is $U W_{\la} = \wt{W}_{\la}
U$) and the set of operators $X\in B(\Hil; \wt{\Hil})$ such that for all  $n \in \Nr$
\[ \sum_{\la \in \La^n} V_{\la} X \wt{V}_{\la}^* =X.\]
\end{cor}

\section{Applications to non-selfadjoint higher-rank graph operator algebras and states on higher-rank graph algebras}

In this last section we present several examples of applications of the Poisson-type transforms to the analysis of the structure
of the related graph operator algebras. In particular we discuss character spaces of the Hardy-type algebras related to a
higher-rank graph $\La$ and characterise purity of a state on $\Cun$ in terms of the related families of unital completely
positive maps.

\subsection*{Character spaces of $\Har$ and $\wHar$}

Let $\cont$ denote the set of all $\La$-contractions on the one dimensional Hilbert space $\bc$.
Note that if $\Vcont \in \cont$, then the condition (iv) in Definition \ref{Lcont} implies that there exists
$a \in \La^0$ such that $V_{a}=1$, $V_b=0$ for $b \in \La^0, b \neq a$. Moreover condition (ii) forces
$V_{\la} =0 $ unless $r(\la)=s(\la) = a $. This together with the remarks after Definition \ref{Lisom}
means that each element of $\cont$ can be identified with
a tuple: $(a, (\alpha_{\la}^{(1)})_{\la \in \ind_1}, \ldots, (\alpha^{(r)}_{\la})_{\la \in \ind_r})$, where
$a\in \La^0$,
for $j \in \oner$ the sequence of complex numbers $(\alpha^{(j)}_{\la})_{\ind_j}$ (where the index set
$\ind_j:=\La^{e_j}_{a,a}$) satisfies $\sum_{\la \in \ind_j} |\alpha^{(j)}_{\la}|^2 \leq 1$,
and for $i, j \in \oner, i \neq j$, there are commutation relations between $\alpha^{(j)}_{\la}$ and
$\alpha^{(i)}_{\mu}$ enforced by the structure of $\La$. To understand the situation  it is essentially sufficient to
determine what happens
if $\La$ has only one vertex (as we are concerned with each $a$ separately). If additionally $r=2$ and $\La$ is finite
we are exactly in the framework considered in \cite{Steve}.

It is easy to see from the above discussion that $\cont$ equipped with the topology of pointwise convergence (that is a net of
$\La$-contractions $\Vcont^{(i)}$ converges to a $\La$-contraction $\Vcont$ if and only if $V^{(i)}_{\la}$ converges to $V_{\la}$
for each $\la \in \La$) is a compact Hausdorff space. Due to the identifications above, $\cont$ can be actually homeomorphically
embedded into the disjoint union of Hilbert spaces $l^2(\ind_j)$, $j \in \oner$, with each factor equipped with the weak
topology. In particular if $\La$ is finite, $\cont$ may be identified with a subset of $\bc^n$ for some $n \in \bn$. The analysis
of the resulting set has been crucial for the full classification of operator algebras $\Har$ associated with a rank-2 graph
having one vertex and finitely many edges of each colour carried out in \cite{SteveKribs} and \cite{Steve}.

\begin{tw} \label{charspace}
Suppose that $\La^0$ is finite. Then the character space of the Banach algebra $\Har$ is homeomorphic to $\cont$.
\end{tw}
\begin{proof}
Suppose first that $f:\Har\to \bc$ is a character and put $V_{\lambda}^f = f(L_{\lambda})$ ($\lambda \in \La$). We claim that the
family $\Vcont^f$ (where complex numbers are viewed as operators on $\bc$) is a $\La$-contraction. The only condition that has to
be checked is whether  $\sum_{\la \in \La^n} |V_{\la}^f|^2 \leq 1$ for all $n \in \Nr$. This is easy to see if we remember that
every character is completely contractive. Let $k \in \bn$ and $\la_1, \cdots, \la_k \in \La^n$. Then
\[ \left(\sum_{i=1}^k |V_{\lambda_i}^f|^2\right)^{\frac{1}{2}}
=  \left\| \left(\begin{bmatrix} V^f_{\lambda_1} & 0 & \cdots  \\
                        \vdots          &    \vdots  & 0 \\
                        V^f_{\lambda_k} & 0 & \cdots
                        \end{bmatrix}\right)\right\| =
                        \left\| f^{(k)}
                        \left(\begin{bmatrix} L_{\lambda_1} &  0 & \cdots  \\
                        \vdots          &    \vdots  & \vdots \\
                        L_{\lambda_k} & 0 & \cdots
                        \end{bmatrix}\right) \right\| \]
\[
\leq \left\|\begin{bmatrix} L_{\lambda_1} &  0 & \cdots  \\
                        \vdots          &    \vdots  & \vdots \\
                        L_{\lambda_k} & 0 & \cdots
                        \end{bmatrix} \right\| =
                        \left\|\sum_{i=1}^k L_{\lambda_i} L_{\lambda_i}^* \right\|^{\frac{1}{2}} \leq 1. \]

Given a $\La$-contraction $\Vcont \in \cont$, define a functional on $\Lin\{L_{\lambda}: \lambda \in \La\}$ by
$f_{\Vcont} (L_{\lambda}) = V_{\lambda}$. Corollary \ref{vNi} implies that $f_{\Vcont}$ is bounded by 1, and the last
statement in Theorem \ref{tPoisson} (together with Remark \ref{notal}) shows that its continuous extension to
$\Har$ is a character.

The correspondences above are in an obvious way inverses of each other
(so that $\Vcont^{f_{\Vcont}} = \Vcont$, $ f_{\Vcont^f} = f$). It is easy to see that the one mapping the character
space into $\cont$ is continuous. By compactness of the sets in question it has to be the desired homeomorphism.
\end{proof}

When $\La^0$ is infinite, it may happen that $f(L_a)=0$ for all $a\in \La^0$ and yet $f\in \Har^*$, $f(1)=1$.
The following is an easy observation which can be proved using the same methods as for the theorem above.

\begin{cor} \label{charinf}
Suppose that $\La^0$ is infinite. The character space of the Banach algebra $\Har$ is homeomorphic to the disjoint topological
union $\cont \cup \{f_0\}$, where $f_0$ is a point representing the character given by $f_0(1)=1$, $f_0 (L_{\la})= 0$ for $\la
\in \La$.
\end{cor}

In the next theorem we characterise the space of weakly-operator (wo-)continuous characters on $\wHar$.

\begin{tw}
There is a one-to-one correspondence between the set of wo-continuous characters on $\wHar$ and those $\La$-contractions
$\Vcont$ in $\cont$ for which
\begin{equation}\label{non1} \kappa_j:=\sum_{\la \in \La^{e_j}} |V_{\la}|^2 < 1\end{equation}
 for each $j \in \oner$.
\end{tw}
\begin{proof}
Note first that the class of wo-continuous characters on $\wHar$, denoted further by ${\cont}_w$,
is a subclass of the family of all continuous characters on $\Har$. If $\La^0$ is infinite, the special character
$f_0$ is not wo-continuous, so that it remains to check which of the elements in $\cont$ determine characters in
${\cont}_w$. To this end we will apply the method from \cite{SteveBaruch}.

Suppose first that $\Vcont \in \cont$ satisfies \eqref{non1}. Let $a \in \La^0$ be such that $V_a=1$.
As for each $n \in \Nr$
\begin{align*} \sum_{\la \in \La^n_{a,a}} |V_{\la}|^2 &=
\sum_{\la_1 \in \La^{n_1 e_1}_{a,a}} \cdots  \sum_{\la_r \in \La^{n_r e_r}_{a,a}} |V_{\la_1} \cdots
V_{\la_r}|^2 \\&=  \left( \sum_{\la_1 \in \La^{e_1}_{a,a}} |V_{\la_1}|^2 \right)^{n_1}
\cdots \left( \sum_{\la_r \in \La^{e_r}_{a,a}} |V_{\la_r}|^2 \right)^{n_r} =
\kappa_1^{n_1} \cdots \kappa_r^{n_r},\end{align*}
we have
\[ \sum_{\la \in \La_{a,a}} |V_{\la}|^2 = \sum_{n\in \Nr} \sum_{\la \in \La^n_{a,a}} |V_{\la}|^2=
\sum_{n\in \Nr} \kappa_1^{n_1} \cdots \kappa_r^{n_r} = \prod_{i=1}^r (1- \kappa_i)^{-1}.\]
Define $\xi \in l^2(\La)$ by
\[ \xi = \sum_{\la \in \La_{a,a}} \overline{V_{\la}} \delta_{\la}\]
and let $\wt{\xi} = \frac{\xi}{\|\xi\|}$.
It is now easy to check that, in the notation of Theorem \ref{charspace},
\begin{equation} \label{vectst} f_{\Vcont} = \langle \wt{\xi}, \cdot \, \wt{\xi}\ra.\end{equation}
Indeed, if $\mu \in \La \setminus \La_{a,a}$ then $L_{\mu} \xi = 0$ or $L_{\mu}^* \xi=0$ and
therefore \[ \langle \wt{\xi}, L_{\mu} \wt{\xi}\ra =0 = f_{\Vcont} (L_{\mu}).\]
If $\mu \in \La_{a,a}$ then
\begin{align*} \langle \xi, L_{\mu} \xi\ra &= \sum_{\la \in \La_{a,a}} \sum_{\gamma \in \La_{a,a}}
\langle \ol{V}_{\lambda} \delta_{\la}, \ol{V}_{\gamma} L_{\mu} \delta_{\gamma} \rangle =
\sum_{\alpha \in \La_{a,a}} \sum_{\gamma \in \La_{a,a}} \langle
 \ol{V}_{\mu \alpha} \delta_{\mu \alpha}, \ol{V}_{\gamma} \delta_{\mu \gamma} \rangle \\&=
\sum_{\alpha \in \La_{a,a}} V_{\mu} |V_{\alpha}|^2 = V_{\mu} \|\xi\|^2 = f(L_{\mu}) \|\xi\|^2.\end{align*} This shows that the
formula \eqref{vectst} holds and therefore $f_{\Vcont}\in {\cont}_w$.

Suppose now that $\Vcont \in \cont$ and $f_{\Vcont} \in {\cont}_w$ and let again
$a \in \La^0$ be such that $V_a =1$. By the arguments similar to those of
\cite{SteveBaruch} we can show that for each $j\in \oner$ restrictions of $f_{\Vcont}$ to the wo-closed algebras
generated by $I$ and $\{L_{\la}: \la \in \La_{a,a}^{e_j}\}$ yield wo-continuous characters on isomorphic copies
of the non-commutative analytic Toeplitz algebras of the type considered in \cite{DavPit}. Theorem 2.3 of that paper
implies that $\Vcont$ has to satisfy \eqref{non1}.
\end{proof}

The theorems above may be used to establish the fact that many non-selfadjoint (higher rank) graph algebras
are not isomorphic (\cite{Pderiv}, \cite{SteveKribs}).

\subsection*{Amenability of $\Har$}

To analyse amenability of the operator algebra $\Har$ we use the characterisation of bounded derivations with values in $\bc$
equipped with a certain natural $\Har$-bimodule structure, obtained by G.Popescu in \cite{Pderiv} (and exploited also in
\cite{PPois}). The main difference here is that in our context we can consider several bimodule structures on $\bc$.

Let $a,b\in \La^0$. Note that the formulas $V_{a} =1$, $V_{\lambda}=0$ for $\lambda \in \La \setminus \{a\}$
define a $\Lambda$ contraction $\Vcont \in \cont$. We will denote the character $f_{\Vcont}$
(see the proof of Theorem \ref{charspace}) simply by $f_a$.
 Let $\bc_{a,b}$ denote $\bc$ equipped with the following actions of $\Har$:
\[ x \cdot \lambda = f_a(x) \lambda, \;\;\; \lambda \cdot x = \lambda f_b(x),  \;\;\; x \in \Har, \lambda \in \bc.\]
It is immediate to check that $\bc_{a,b}$ equipped with these actions is a (dual, normal) Banach bimodule.

Denote $\La^{1}_{a,b} = \{ \la \in \Lambda: |\sigma (\lambda)| =1, r(\lambda)=a, s(\lambda)=b\}$.

\begin{lem}
Let $a,b \in \La^0$, $a\neq b$. Every derivation $\delta: \Har \to \bc_{a,b}$ is given by a linear extension of
the formula
\[ \delta (L_{\lambda}) = \alpha_{\lambda}, \;\;\; \lambda\in \La^{1}_{a,b},\]
\[ \delta (L_b) = - \delta(L_a) = \alpha_b,\]
\[ \delta(1)=\delta (L_{\lambda}) = 0, \;\;\; \lambda\in \La \setminus \{\La^{1}_{a,b}\cup \{a,b\}\},\]
where $\alpha_a \in \bc$, and the family  $\{\alpha_{\lambda}\in \bc:\lambda\in \La^{1}_{a,b}\}$ is such that
 $\sum_{\lambda \in \La^{1}_{a,b}} |\alpha_{\lambda}|^2 < \infty$.
 Conversely, every $\alpha \in \bc$ and square integrable family $\{\alpha_{\lambda}:\lambda\in \La^{1}_{a,b}\}$
 define a derivation via the formulas above.
The same remains true for $a=b$, except that $\alpha_b=0$.
\end{lem}
\begin{proof}
Suppose first that $a \neq b$ and we are given a bounded derivation $\delta: \Har\to \bc_{a,b}$.
Then the formulas
\[ 0 = \delta(0) = \delta (L_a L_b) = \delta(L_a) + \delta(L_b),\]
\[ 0 = \delta(L_c L_b) = \delta(L_c)\, 1 + 0 \, \delta (L_b), \;\;\; c \in \La^0, c \notin \{a,b\},\]
\[ \delta(L_{\mu \nu}) = \delta(L_{\mu}) \,0 + 0 \,\delta(L_{\nu})=0, \;\;\; \mu, \nu \in \La \setminus \La^0,\]
\[ \delta (L_{\mu}) = 0 \,L_{\mu} + L_{r(\mu)}\, 0 = 0, \;\;\;\; \mu \in \La, |\sigma(\mu)|=1, r (\mu) \neq a,\]
\[ \delta (L_{\mu}) = L_{\mu} \,0 + 0 \,L_{s(\mu)}  = 0, \;\;\;\; \mu \in \La, |\sigma(\mu)|=1, s (\mu) \neq b,\]
imply that $\delta (L_{\la}) = 0$ unless $\la\in \La^{1}_{a,b} \cup \{a,b\}$.
Put $\alpha_b=\delta(L_b)$, $\alpha_{\la}= \delta (L_{\la})$ for $\la\in \La^{1}_{a,b}$.
Choose $j\in \oner$ and consider the set $\La^{e_j}_{a,b}:=\{\la \in \La^{e_j}: r(\la) = a, s(\la)=b\}$.
As $\delta$ is bounded we have in particular for each finite set $F \subset \La^{e_j}_{a,b}$ and all coefficients
$\gamma_{\la} \in \bc$ ($\la \in F$)
\[ \left| \sum_{\la \in F}  \gamma_{\la} \alpha_{\la}\right|
= \left|\delta \Big(\sum_{\la \in F} \gamma_{\la} L_{\la}\Big)\right|
\leq \|\delta\| \|\sum_{\la \in F} \gamma_{\la} L_{\la}\|.\]
It remains to compute the norm on the right hand side. This however is easy if we note that due to the
factorisation property the operator
$\sum_{\la \in F} \gamma_{\la} L_{\la}$ has a certain block structure with respect to the orthogonal decomposition
$l^2(\La) = \bigoplus_{n \in \Nr} \La^n$. More precisely, for any $\xi \in l^2(\La)$,
$\xi = \sum_{\mu \in \La} \beta_{\mu} \delta_{\mu}$ we have
\begin{align*}
 \left\|\sum_{\la \in F} \gamma_{\la} L_{\la} (\xi) \right\|^2
&= \left\|\sum_{\la \in F} \gamma_{\la} L_{\la} (\sum_{\mu \in \La} \beta_{\mu} \delta_{\mu}) \right\|^2 =
  \left\| \sum_{\la \in F, \mu \in \La}\gamma_{\la} \beta_{\mu} \delta_{\la \mu} \right\|^2
  \\&=
  \sum_{\la \in F, \mu \in \La} |\gamma_{\la} \beta_{\mu}|^2 =
  \sum_{\la \in F} |\gamma_{\la}|^2 \sum_{\mu \in \La} |\beta_{\mu}|^2 =
\sum_{\la \in F} |\gamma_{\la}|^2 \|\xi\|^2.\end{align*}
This implies that
\[ \left| \sum_{\la \in F}  \gamma_{\la} \alpha_{\la}\right|
\leq \|\delta\| \left(\sum_{\la \in F} |\gamma_{\la}|^2\right)^{\frac{1}{2}}.\]
But this means that $\sum_{\la \in F} |\alpha_{\la}|^2 \leq \| \delta\|^2$, and as $F$ was an arbitrary finite subset
of $\La^{e_j}_{a,b}$, and $\La^1_{a,b} = \bigcup_{j=1}^r \La^{e_j}_{a,b}$ we proved that
$\sum_{\la \in \La^1_{a,b}} |\alpha_{\la}|^2 \leq r\| \delta\|^2$.

The converse implication is even easier to establish. Suppose that we are given $\alpha_b \in \bc$ and a square summable family
$\{\alpha_{\lambda}:\lambda\in \La^{1}_{a,b}\}$. Define a map $\delta$ by the (linear extension of the) formulas in the lemma. As
it is easy to check that $\delta$ is a derivation, it remains to show that it is bounded. To this end let $\{\gamma_{\la}: \la
\in \La\}$ be a finitely supported family of complex numbers and consider $f = \sum_{\la \in \La} \gamma_\la L_{\la} \in \Har$.
As $\|f (\delta_{b})\|= \|\sum_{\la \in \La: s(\la)=b} \gamma_{\la} \delta_{\la}\| = \left(\sum_{\la \in \La: s(\la)=b}
|\gamma_{\la}|^2\right)^{\frac{1}{2}}$ we have
\[\| f \|\geq \left(\sum_{\la \in \La: s(\la)=b} |\gamma_{\la}|^2\right)^{\frac{1}{2}}.\]
Analogously \[\| f \|\geq \left(\sum_{\la \in \La: s(\la)=a} |\gamma_{\la}|^2\right)^{\frac{1}{2}}.\]
This allows to obtain the following estimate:
\begin{align*}|\delta (f)|  &=
\left|\gamma_b \alpha_b - \gamma_a \alpha_b + \sum_{\la \in \La^1_{a,b}} \gamma_{\la} \alpha_{\la}\right| \leq |\gamma_b
\alpha_b| + \left(\sum_{\la \in \La^1_{a,b} \cup \{a\}} |\alpha_{\la}|^2 \right)^{\frac{1}{2}} \left(\sum_{\la \in \La^1_{a,b}
\cup \{a\}} |\gamma_{\la}|^2 \right)^{\frac{1}{2}} \\& \leq \left( |\alpha_b| + \left(\sum_{\la \in \La^1_{a,b} \cup \{a\}}
|\alpha_{\la}|^2 \right)^{\frac{1}{2}} \right) \|f\|.
\end{align*}
This ends the proof in the case $a \neq b$.

It is easy to see that if $a=b$, then $\alpha_b$ above has to be equal $0$, and all the other arguments need not be changed.

\end{proof}

\begin{cor}
Let $a,b \in \La^0$. The (first continuous) cohomology group $H_c^1(\Har, \bc_{a,b})$ is isomorphic to the space
$l^2(\La^1_{a,b})$.
\end{cor}
\begin{proof}
The corollary follows from the lemma above and the fact that all inner derivations $\delta:\Har \to \bc_{a,a}$ are
trivial, and for $a \neq b$ all inner derivations $\delta:\Har  \to \bc_{a,b}$ are given by the formula
\[ \delta( L_a) = - \delta(L_b) = \alpha,\]
\[ \delta (L_{\lambda}) = 0, \;\;\; \lambda\in \La \setminus \{a, b\},\]
where $\alpha \in \bc$.
\end{proof}

This gives the following result:

\begin{tw}
Suppose that $\Lambda$ is  finitely aligned  and nontrivial (i.e.\ there exists $\lambda\in \La$ such that
$\sigma(\lambda) \neq 0$). Then the algebra $\Har$ is not amenable.
\end{tw}
\begin{proof}
Follows immediately from the lemma above.
\end{proof}

Note that this provides a wide family of examples of both finite- and infinite-dimensional non-amenable algebras (the
simplest being given by a subalgebra of $M_3(\bc)$ constructed from a rank-one graph with two vertices and a
connecting edge).

\subsection*{States on $\Cun$}

The first result is a natural generalisation of Theorem 2.1 of \cite{pureod}. For its formulation we introduce a global neutral element $\emptyset$ for $\La$ verifying $\la \emptyset = \la$ for all $\la \in \La$.

\begin{tw}  \label{state}
Let $\La$ be finite and with no sources. There is a 1-1 correspondence between the following objects:
\begin{rlist}
\item states $\omega:\Cun \to \bc$;
\item positive definite kernels $k:(\La \cup \{\emptyset\}) \times (\La \cup
\{\emptyset\}) \to \bc$ such that $k (\emptyset, \emptyset) =1$ and for all $\mu,\nu \in \La \cup \{\emptyset\}$, $n \in \Nr$
\[ \sum_{\la \in \La^n} k(\mu\la, \nu \la) = k(\mu,\nu);\]
\item unitary equivalence classes of the triples $(\Kil, \Omega,
\Vcont)$, where $\Kil$ is a Hilbert space, $\Omega \in \Kil$ has norm 1, $\Vcont$ is a $\La$-isometry on $\Kil$
and $\Kil = \overline{\textup{Lin}}\{V_{\la}^* \Omega: \la \in \La\}$.
\end{rlist}

The correspondence is given by the formulas
\begin{equation} \label{corr} \omega(s_{\la} s_{\mu}^*) = k(\la, \mu) = \langle V_{\la}^*
\Omega,V_{\mu}^* \Omega \ra, \;\;\; (V_{\emptyset}:= I_{\Kil}, s_{\emptyset}:=1_{\Cun}) \end{equation}
\end{tw}

\begin{proof}
It is clear that whenever a state $\omega$ on $\Cun$ and a triple $(\Kil, \Omega,
\Vcont)$ satisfying the conditions in (iii) is given, the function
$k:(\La \cup \{\emptyset\}) \times (\La \cup
\{\emptyset\}) \to \bc$ defined via the formulas in \eqref{corr} is a positive definite kernel.
Given a positive definite kernel $k$ as in (ii), the standard Kolmogorov construction
(see for example \cite{Kolm}) yields a Hilbert space $\Kil$ (unique up to an isomorphism)
and a map $T:(\La \cup \{\emptyset\})  \to \Kil$ such that
\[ \langle T(\la), T(\mu) \ra = k(\la, \mu), \;\;\; \la, \mu \in \La \cup \{\emptyset\},\]
and $\Kil = \ol{\Lin}\{T(\la): \la \in  \La \cup \{\emptyset\}\}$.
Define for each $\mu \in \La$
\[S(\mu) (T(\la)) = T( \la \mu),  \;\;\; \la \in \La \cup \{\emptyset\} \]
(with the convention that if $\la \in \La$ and $s(\mu)\neq r(\la)$
then $T(\mu \la) = 0$). The linear extension of $S(\mu)$ is contractive: let $n = \sigma(\mu)$,
$k \in \bn$, $\la_1, \ldots, \la_k \in \La \cup \{\emptyset\}$,
$\alpha_1, \ldots, \alpha_k \in \bc$, and compute:
\begin{align*} \Big\| S_{\mu} \Big( &\sum_{i=1}^k \alpha_i T(\la_i) \Big) \Big\|^2 =
\Big\| \sum_{i=1}^k \alpha_i T( \la_i \mu) \Big\|^2
\\&= \sum_{i,j=1}^k \ol{\alpha}_i \alpha_j \langle  T(\la_i\mu), T( \la_j\mu) \rangle =
\sum_{i,j=1}^k \ol{\alpha}_i \alpha_j   k( \la_i \mu,  \la_j\mu) \\& \leq
\sum_{i,j=1}^k \sum_{\gamma \in \La^n} \ol{\alpha}_{i} \alpha_j   k(\la_i \gamma, \la_j \gamma)
= \sum_{i,j=1}^k \ol{\alpha}_i \alpha_j   k(\la_i, \la_j) =
\Big\| \sum_{i=1}^k \alpha_i T(\la_i) \Big\|^2
\end{align*}
(note that in the inequality above we used the fact that $k$ is a positive-definite kernel).
Put $V_{\mu} = S(\mu)^*$ and $\Omega= T(\emptyset)$. Then $\Vcont$ is a $\La$-isometry and
the triple $(\Kil, \Omega, \Vcont)$ satisfies the conditions in (iii) and is compatible with
$k$ (satisfies one of the equalities in \eqref{corr}).

It remains to show that every triple $(\Kil, \Omega, \Vcont)$ in (ii) yields a state on $\Cun$.
This however is immediate via the Poisson transform from Theorem \ref{cPoisson}:
define $\omega$ by the formula
\[ \omega(X) = \langle \Omega, \cPois(X) \Omega\ra, \;\;\; X \in \Cun.\]
It is easy to check that the compatibility conditions in \eqref{corr} hold.
\end{proof}

As noted in \cite{pureod}, there is a direct way of constructing the triple  $(\Kil, \Omega,
\Vcont)$ corresponding to a given state $\omega$ on $\Cun$. It is described in the next remark.

\begin{rem}    \label{canon}
Let $\La$ be finite and with no sources. Suppose that $\omega$ is a state on $\Cun$.
Let $(\pi,\Hil_{\omega}, \Omega)$ be the corresponding  GNS triple and define the $\La$-isometry
$\Wcont$ on $\Hil_{\omega}$ by putting $W_{\la} = \pi(s_{\la})$ ($\la \in \La$). Let
$\Kil = \ol{\Lin}\{W_{\la}^* \Omega : \la \in\La \}$ and let $P$ denote the orthogonal projection
from $\Hil_{\omega}$ onto $\Kil$. Put $V_{\la} = P W_{\la} P$ for $\la \in \La$.
Then it is easy to check that the resulting family $\Vcont$ is a $\La$-isometry on $\Kil$. For example
if $\la, \mu \in \La$ then for arbitrary $ \nu, \gamma \in \La$
\[ \langle W_{\nu}^* \Omega, V_{\la} V_{\mu} W_{\gamma}^* \Omega \ra =
\langle  V_{\mu}^* W_{\nu \la }^* \Omega,  W_{\gamma}^* \Omega \ra =
\langle ( W_{\la} W_{\mu})^* W_{\nu \la }^* \Omega,  W_{\gamma}^* \Omega \ra.\]
Now the latter is equal $0$ if $r(\mu) \neq s(\la)$ and if $r (\mu) =s(\la)$
\[ \langle W_{\nu}^* \Omega, V_{\la} V_{\mu} W_{\gamma}^* \Omega \ra
=\langle ( W_{\la \mu})^* W_{\nu \la }^* \Omega,  W_{\gamma}^* \Omega \ra =
\langle   W_{\nu \la }^* \Omega,  V_{\la \mu} W_{\gamma}^* \Omega \ra,\]
which implies that $V_{\la} V_{\mu} = V_{\la \mu}$ if  $r(\mu) = s(\la)$
and $V_{\la} V_{\mu} = 0$ otherwise. In particular this implies that each $V_a$ ($a \in \La^0$) is
an idempotent (as it is obviously a contraction, it has to be a selfadjoint projection).
Moreover $\Wcont$ is the minimal
dilation of $\Vcont$ to a Cuntz-Pimsner $\La$-family. This follows immediately from
the fact that each element in $\Hil_{\omega}$ can be approximated by linear combinations
of vectors of the type $W_{\la} W_{\mu}^* \Omega$ ($\la, \mu \in \La$) and the latter are obviously
linear combinations of vectors of the form $W_{\la} \xi$ ($\la\in \La, \xi \in \Kil$).
\end{rem}

The following lemma is an immediate consequence of Theorem 3.9 and the discussion above.
We use the notation introduced in Definition \ref{defcp}.

\begin{lem}      \label{fix2}
Suppose that $\La$ is finite and has no sources. Let $\omega$ be a state on $\Cun$ and let
$(\pi,\Hil_{\omega}, \Omega)$, $\Kil$, $P$,
$\Wcont$ and $\Vcont$ be defined as in Remark \ref{canon}.
Then the map $X \longrightarrow PXP$ yields a completely isometric complete order isomorphism
between $\pi(\Cun)'$ and the  operator system
$\Fix \sigma_{\Vcont}$.
\end{lem}

The following theorem extends the characterisation of pure states on the Cuntz algebras ${\mathcal O}_d$ given in \cite{pureod}. The states are characterised in terms of the properties
of a corresponding $\La$-isometry.

\begin{tw}
Suppose that $\La$ is finite and has no sources. Let $\omega:\Cun \to \bc$ be a state, let $(\pi,\Hil, \Omega)$ be the GNS triple
for $\omega$, let $\Wcont$ be the $\La$-isometry  given by $W_{\la} = \pi (s_{\la}) (\la \in \La)$. Denote the closed subspace
spanned by $W_{\la}^* \Omega$ ($\la \in \La$) by $\Kil$, let $P$ denote the orthogonal projection from $\Hil$ to $\Kil$ and let
$\Vcont$ be the $\La$-isometry on $\Kil$ constructed according to Remark \ref{canon}. The following are equivalent:
\begin{rlist}
\item $\omega$ is pure;
\item $ \Fix \sigma_{\Wcont} = \bc I_{\Hil}$;
\item $ \Fix \sigma_{\Vcont}  = \bc I_{\Kil}$;
\item $\Vcont$ acts irreducibly on $\Kil$ and $P \in \pi(\Cun)''$.
\end{rlist}
\end{tw}
\begin{proof}
%Recall first that Lemma
%\ref{same} implies that $\Fizx \sigma_{\Wcont} = \Fix \sigma_{\Wcont}$,
%$ \Fizx \sigma_{\Vcont} = \Fix \sigma_{\Vcont}$.

The state $\omega$ is pure if and only if  $\pi (\Cun)' =\bc I_{\Hil}$. By Lemma \ref{fix1}
the latter is equivalent to $\Fix \sigma_{\Wcont} = \bc I_{\Hil}$.
Similarly  by Lemma \ref{fix2}  $\pi (\Cun)' =\bc I_{\Hil}$ if and only if
the selfadjoint part of $\Fix \sigma_{\Vcont}$ is one-dimensional (over reals).
This happens if and only if $\Fix \sigma_{\Vcont} = \bc I_{\Kil}$ and  the equivalences
(i)$\Leftrightarrow$(ii)$\Leftrightarrow$(iii) have been proved.  If (i) holds then
$\pi (\Cun)'' = B(\Hil)$ and the commutant of $\Vcont$ is contained
in $\Fix \sigma_{\Vcont}$ so by (iii) has to be trivial. Thus (i) and (iii) imply (iv).

Suppose finally that (iv) holds.
Since $P \in \pi(\Cun)''$ it is easy to see that
$P \pi (\Cun)' P $ is an algebra. By Lemma \ref{fix2}
 $\Fix \sigma_{\Vcont} = P \pi (\Cun)' P $.
Suppose that $p \in B(\Kil)$ is a projection.
Note that as $\Fix \sigma_{\Vcont} = \bigcap_{j=1}^r \Fix \sigma_{\Vcont} (e_j) \cap
\Fix \sigma_{\Vcont} (0)$ and
the commutant of $\Vcont$ is the intersection of  the commutants of the sets
$\{V_{\la}, V_{\la}^*: \la \in \La^{e_j}\}$   ($j \in \oner$) and the commutant of
$\{V_a: a \in \La^{0}\}$, Lemma 3.3 of \cite{pureod} implies that $p$ belongs to
 $\Fix \sigma_{\Vcont}$
if and only if it belongs to the commutant of $\Vcont$. This implies that the largest
von Neumann subalgebra (or, equivalently here, the largest $*$-subalgebra) of
$\Fix \sigma_{\Vcont}$ coincides with the commutant of $\Vcont$. As we established before
that $\Fix \sigma_{\Vcont}$ is an algebra, (iii) follows.
\end{proof}

In Section 7 of \cite{pureod} a version of Theorem \ref{state} was used to construct translationally
invariant states on two-sided quantum spin chains. The construction can be generalised to our framework. Given
a triple $(\Kil, \Omega, \Vcont)$ as in Theorem \ref{state} (iii) we obtain a state on $\Cun$ and thus also
a state on the canonical AF algebra  $\Gauge$ (assume $\La$ is finite) arising as a fixed point subalgebra for the gauge action
on $\Cun$ (see \cite{kupa} for the details). The resulting state can be interpreted as being translationally invariant
with respect to the shifts in each of the $r$ available directions. Note however that already for basic rank-2 examples
the situation becomes very subtle. When $\La$ is rank-2 graph which has one vertex, two `red' edges $e_1, e_2$
and two `green' edges $f_1, f_2$ with the factorisation rules $e_1 f_2 = f_1 e_2$, $e_2 f_1 = f_2 e_1$,
the algebra $\Cun$ is isomorphic to $C(\bt) \ot \mathcal{O}_2$, and $\Gauge \approx UHF(2^{\infty})$. It is easy to see
that the states on $UHF(2^{\infty})$ arising via the construction described above
are quite different to those produced in \cite{pureod}.

\subsection*{ACKNOWLEDGMENT}
The work on this paper was started during the first-named author's visit to the Indian Statistical Institute in Bangalore in December 2006/January 2007. AS would like to thank Professor B.V.\,Rajarama Bhat for the invitation and providing ideal research conditions.

\end{document}